\documentclass{siamltex}
\usepackage{graphicx}
\usepackage{color}
\usepackage{subfigure}
\usepackage{captcont}
\usepackage{subfloat}

\usepackage{amsmath}
\usepackage{amssymb}
\usepackage{amsbsy}
\usepackage{comment}

\title{An adaptive, high-order phase-space remapping for the two-dimensional Vlasov-Poisson equations}

\author{Bei Wang \thanks{Department of Applied Science, University of California, 1 Shields Avenue, Davis CA 95616.
Present Address: Princeton Institute for Computational Science and Engineering, Princeton University, Princeton, NJ 08540. E-mail: beiwang@princeton.edu} 
\and
Greg Miller\thanks{Department of Chemical Engineering and Materials Science, 1 Shields Avenue, Davis CA 95616.}
\and
Phil Colella\thanks{Applied Numerical Algorithms Group, Lawrence Berkeley National Laboratory, MS 50A-1148, 1 Cyclotron Road, Berkeley, CA 94720.}}

\begin{document}
\maketitle


\begin{abstract}
The numerical solution of high dimensional Vlasov equation is usually performed by particle-in-cell (PIC) methods. 
However, due to the well-known numerical noise, it is challenging to use PIC methods to get a precise description of the distribution function in phase space.
To control the numerical error, we introduce an adaptive phase-space remapping which regularizes the particle distribution by periodically reconstructing the distribution function on a hierarchy of 
phase-space grids with high-order interpolations. The positivity of the distribution function can be preserved by using a local 
redistribution technique. The method has been successfully applied to a set of classical plasma problems in one dimension. 
In this paper, we present the algorithm for the two dimensional Vlasov-Poisson equations. An efficient Poisson solver with 
infinite domain boundary conditions is used. The parallel scalability of the algorithm on massively parallel computers will be discussed.
\end{abstract}



\pagestyle{myheadings}
\thispagestyle{plain}
\sloppy 

\section{Introduction}
The Vlasov equation describes the dynamics of a species of charged particles under electromagnetic fields. In the electrostatic case, the normalized equation reads
\begin{equation}
 \frac{\partial f}{\partial t} + \boldsymbol{v} \cdot \nabla_{\boldsymbol{x}} f +
{(-1)}^s (\boldsymbol{E} + \boldsymbol{E}^e) \cdot \nabla_{\boldsymbol{v}} f = 0,
\label{eqn:vlasov}
\end{equation}
where $f(\boldsymbol{x}, \boldsymbol{v}, t)$ is the distribution function of the species in phase space 
$(\boldsymbol{x}, \boldsymbol{v}) \in \mathbb{R}^d\times\mathbb{R}^d$ with $d={1,2,3}$. $s$ is $0$ for positive charges and is $1$ otherwise. $\boldsymbol{E}$ and 
$\boldsymbol{E}^e$ denote the self-consistent and the external electric field, respectively.
This equation is the simplest model to study collisionless plasmas and beam propagation which is of importance to controlled thermonuclear fusion and accelerator modeling. 

The Vlasov equation is a nonlinear hyperbolic equation in phase space so methods of solution can be guided by the 
well-established numerical analysis of classical partial differential equations. Accordingly, grid methods in fluid dynamics, such as transform methods, 
finite-volume methods, and semi-Lagrangian methods, can be employed. Operator splitting was successfully applied to the solution of the Vlasov equation
by Cheng and Knorr \cite{Cheng76} in 1970s. It reduces the solution of the multi-dimensional Vlasov equation to a set of 
one-dimensional advection problems, and therefore has become a widely used technique. Even with these well-established algorithms, performing high-dimensional 
simulations using grid methods is still a challenging task. The issue is the computational time and memory cost in dealing with the whole six dimensional 
phase space. With the advances of supercomputer, grids methods have achieved large development in the last decade.
In semi-Lagrangian methods, Sonnendrucker et al. \cite{Sonnendrucker98} introduced the cubic spline method. Nakamura and Yabe \cite{Nakamura99} 
introduced the cubic interpolated propagation method. In finite-volume methods, Fijalkow \cite{Fijalkow99} presented 
the flux balance method. A similar idea is used in a high-order
finite-volume method based on mapped coordinate by Colella, Dorr and Hittinger \cite{Colella10}. Filbet, Sonnendrucker, and Bertrand \cite{Filbet01} proposed the positive and flux conservative scheme 
using the idea of limiter.

A more widely used approach for the solution of the Vlasov equation is PIC methods \cite{Hockney81,Birdsall91}. In PIC methods, 
the particles, a Lagrangian discretization of the distribution function, follow trajectories computed from the characteristic curves given by the Vlasov equation, whereas the
 self-consistent fields are calculated on a grid. Since the methods employ the fundamental equations without much
 approximation, it allows us to observe most of the physics in a plasma system with relatively few particles.
However, as with all other particle methods, PIC methods suffer from numerical noise such that they have difficulty in simulating some 
problems, e.g., the problem with large dynamic ranges in velocity space. To remedy this deficiency, there are usually two approaches.
One method is the so-called $\delta f$ method \cite{Dimits93,Parker93,Lee11}, discretizing only the perturbation $\delta f$ with respect
to an equilibrium state $f_0$ based on a particle method. The $\delta f$ method has been successfully used in realistic applications, e.g., microtubulence in 
magnetic confined plasmas \cite{Lin98}. The limitation of this method is that it can only be applied to the problems which are close to
equilibrium. An alternative approach is through periodically reconstructing the distribution function on a grid in phase space. 
Such remapping technique has been used in particle methods in fluid dynamics \cite{Cottet00,Chaniotis02}, i.e., vortex methods and smoothed particle hydrodynamics (SPH), to maintain regularity of 
the particle distribution and thereby improve accuracy, but has much more limited use in PIC methods in plasma physics. It is worth mentioning that early work 
of Denavit \cite{Denavit72} and more recent work of Vadlamani \cite{Vadlamani04} and Yang \cite{Chen07} used the idea of remapping for PIC methods. However, They all used low order interpolation
function which results in a first order method overall.

We studied a high order remapping scheme to PIC methods for the solution of the one-dimensional Vlasov-Poisson equations early \cite{Wang10}. 
Meanwhile, we provided a local redistribution technique such that the positivity of the distribution function could be preserved after high-order remapping.
The initial numerical experiments on a set of classical plasma problems in one dimension are very encouraging. Remapping significantly reduces the numerical noise and 
results in a more consistent second-order convergence rate in the electric field error. We also investigated the effects of integrating mesh refinement
to the uniform remapping. This is motivated by the observation that remapping, a numerical diffusive procedure, tends to create 
a large number of small-strength particles at the low density region of the distribution function. Mesh refinement has the potential to reduce this side effect. 

In this paper, we extend the algorithm to the solution of the two-dimensional Vlasov-Poisson equations. This includes the use of an efficient
Poisson solver with infinite domain boundary conditions for beam problems. High-dimensional simulations are very expensive with respect to memory usage. We perform the simulation on a parallel 
machine using domain decomposition in physical space. A scalable implementation based on domain decomposition in phase space will be discussed. 
We consider two types of numerical tests: plasma problems including linear Landau damping and the two stream instability, and a beam 
problem based on the paraxial model \cite{Filbet06}. 

The rest of the paper is organized as follows. In \S\ref{sec:algorithms}, we first review the classical PIC methods for the Vlasov-Poisson equations. An efficient
algorithm which solves the Poisson equation with infinite boundary conditions is described. 
Then we present the high-order and positive remapping on a hierarchy of locally-refined grids in two dimensions.
\S\ref{sec:implementation} discusses the parallel implementation of the algorithm. We show the numerical results in \S\ref{sec:tests}. 
Conclusion and future research will be given at the end.

\section{Algorithms}
\label{sec:algorithms}
\subsection{PIC methods}
PIC methods are based on the Lagrangian description of the Vlasov equation
\begin{equation}
 \frac{df(\boldsymbol{X}, \boldsymbol{V}, t)}{dt} = 0,
 \label{eqn:dfdt}
\end{equation}
where the characteristics $(\boldsymbol{X}(t), \boldsymbol{V}(t))$ are the solution of the equation of motion:
\begin{equation}
 \frac{d\boldsymbol{X}}{dt} = \boldsymbol{V}(t), 
\\
 \frac{d\boldsymbol{V}}{dt} = {(-1)}^s (\boldsymbol{E}(\boldsymbol{X}, t) + \boldsymbol{E}^e(\boldsymbol{X}, t))
\label{eqn:trajectory}
\end{equation}
with initial conditions $\boldsymbol{X}(t=0)=\boldsymbol{x}$ and $\boldsymbol{V}(t=0)=\boldsymbol{v}$.

In the beginning, the distribution function is approximated by a collection of point particles,
\begin{equation}
 f(\boldsymbol{x},\boldsymbol{v}, t=0) \approx \sum_k q_k \delta(\boldsymbol{x}-\boldsymbol{x}_k) \delta(\boldsymbol{v}-\boldsymbol{v}_k),
\label{eqn:pointparticle}
\end{equation}
where $(\boldsymbol{x}_k, \boldsymbol{v}_k)$ is a initial particle location at the cell center of a grid in phase space (quite start). 
$q_{k}=f(\boldsymbol{x}_{k}, \boldsymbol{v}_{k}, t=0) h_xh_yh_{v_x}h_{v_y}$ is the weight of a particle.
Then each particle follows a trajectory described by the equation of motion, 
\begin{equation}
 \frac{dq_k}{dt}=0, \quad
 \frac{d\tilde{\boldsymbol{X}}_k}{dt}(t) = \tilde{\boldsymbol{V}}_k(t), \quad 
 \frac{d\tilde{\boldsymbol{V}}_k}{dt}(t) = {(-)}^s (\tilde{\boldsymbol{E}}_k(t) + \boldsymbol{E}^e_k(t)),
\label{eqn:eom}
\end{equation}
where $\tilde{\boldsymbol{X}}_k(t=0) = \boldsymbol{x}_k$ and $\tilde{\boldsymbol{V}}_k(t=0) = \boldsymbol{v}_k$. 

At any time that a smooth representation of the distribution function is required, we approximate the function with a collection of finite size particles,
where the exact delta function is replaced by a smoothed delta function. That is,
\begin{equation}
 f(\boldsymbol{x},\boldsymbol{v}, t) \approx \sum_k q_k \delta_{\boldsymbol{\varepsilon}_{\boldsymbol{x}}}(\boldsymbol{x}-\tilde{\boldsymbol{X}}_k(t)) 
\delta_{\boldsymbol{\varepsilon}_{\boldsymbol{v}}}(\boldsymbol{v}-\tilde{\boldsymbol{V}}_k(t)), \quad t > 0.
\end{equation}
The smoothed delta function satisfies 
\begin{equation}
\int_{\mathbb{R}^2} \delta_{\boldsymbol{\varepsilon}}(\boldsymbol{y}) d\boldsymbol{y} = 1
\end{equation}
and
\begin{equation}
 \delta_{\boldsymbol{\varepsilon}}(\boldsymbol{y}) = \prod_{d=0}^1 \frac{1}{\varepsilon_d} u\left(\frac{y_d}{\varepsilon_d}\right),
\label{eqn:smoothdelta}
\end{equation}
where $u$ is any interpolation function and $\boldsymbol{\varepsilon}$ is the stencil size. Usually, the stencil size for the smoothed delta function in physical space 
$\boldsymbol{\varepsilon}_{\boldsymbol{x}}$ is chosen as the same as the mesh spacing of the Poisson solver.
The typical interpolation function for PIC methods is the first-order interpolation function
\begin{equation}
 u_1(z) = \begin{cases}
1-|z| & 0 \le |z| \le 1 \\
0 & \mbox{otherwise.} \\
\end{cases}
\label{eqn:linearfunc}
\end{equation}

The flow of a PIC scheme is
\begin{itemize}
 \item Assign particle charges on a grid in physical space,
 \begin{equation}
 \tilde{\rho}(\boldsymbol{x}_{\boldsymbol{j}}, t) = \sum_k q_k \delta_{\boldsymbol{\varepsilon}_{\boldsymbol{x}}}(\boldsymbol{x}_{\boldsymbol{j}} - \tilde{\boldsymbol{X}}_k(t)),
 \label{eqn:rho}
 \end{equation}
where $\boldsymbol{j} \in \mathbb{Z}^2$ are the node index of the Cartesian grid in physical space. 
The grid size is chosen as the same as the stencil size of the smoothed delta function $\boldsymbol{\varepsilon}_{\boldsymbol{x}}$. 
In the case of the first-order interpolation function, for each node $\boldsymbol{j}$, the sum is restricted to the particles with 
$|\boldsymbol{x}_{\boldsymbol{j}} - \tilde{\boldsymbol{X}}_k(t)| \leq \boldsymbol{\varepsilon}_{\boldsymbol{x}}$. 
 \item Solve the Poisson equation on the grid with a second-order finite-difference method:
\begin{equation}
-{(\triangle^H \phi)}_{\boldsymbol{j}} = -\sum_{d=0}^{1} \frac{\phi_{\boldsymbol{j}+\boldsymbol{e}^d} - 2 \phi_{\boldsymbol{j}}
+ \phi_{\boldsymbol{j}-\boldsymbol{e}^d}}{\boldsymbol{\varepsilon}_{\boldsymbol{x}_d}^2} =  {(-1)}^{s} \tilde{\rho}_{\boldsymbol{j}} + {\rho}_{\mbox{background}}
\label{chp1:eqn:discretePoisson}
\end{equation}
and 
\begin{equation}
\tilde{\boldsymbol{E}}_{\boldsymbol{j}}^d = \frac{\phi_{\boldsymbol{j}-\boldsymbol{e}^d} - \phi_{\boldsymbol{j}+\boldsymbol{e}^d}}{2\boldsymbol{\varepsilon}_{\boldsymbol{x}_d}}. 
\end{equation}
$\rho_{\mbox{background}}$ is the background charge density if applicable. With given boundary conditions, the discrete Poisson equation is usually solved by a fast Poisson solver, such as FFTs or multigrid methods.
 \item Interpolate the calculated field back to the particle locations with the same interpolation function in equation (\ref{eqn:smoothdelta}). It is worth mentioning that 
a different interpolation function will introduce self-force errors \cite{ColellaNorgaard10}.
 \item Integrate the equation of motion numerically, for example, using the second-order Runge-Kutta method.
\end{itemize}

\subsection{Solving the Poisson equation with infinite domain boundary conditions}
To model beam problems, the Poisson equation with infinite domain boundary conditions needs to be solved. We compute the solution using a new version of the 
James-Lackner method \cite{James77} by McCorquodale et al. \cite{Mccorquodale05,Mccorquodale07}. This method solves two Dirichlet boundary 
problems plus a boundary to boundary convolution.

\begin{figure}[hbtp]
\begin{center}
\resizebox{12cm}{!}{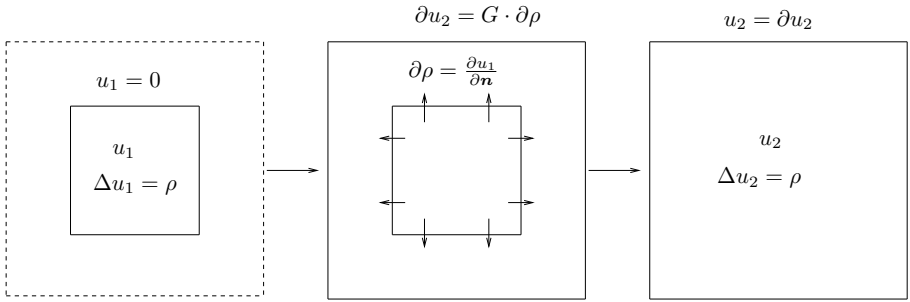}
\caption{James Algorithm \quad Left: Solve for $u_1$ on $D_1$ \quad Middle: Calculate surface charge and perform convolution \quad Right: Solve $u_2$ on $D_2$ 
\label{fig:james} } 
\end{center}
\end{figure}

We briefly describe the algorithm below. Assume $D_0$ is the support domain of the right-hand side $\rho$, we can solve the Poisson equation with infinite domain boundary
conditions on a slightly larger domain $D_1>D$ with inhomogeneous Dirichlet boundary conditions. The boundary value can be calculated by Green's function convolution from
the source $\rho$ to the domain boundary $\partial D_1$. The volume source $\rho$ to boundary $\partial D_1$ convolution is relatively expensive, in particular for 
$3D$ problems. Instead of using a volume to boundary convolution, we can compute the boundary value by performing a boundary $\partial D_1$ to boundary $\partial D_2$ 
convolution and solving another Poisson equation on a domain $D_2>D_1>D$ with Dirichlet boundary condition. The procedure of James' algorithm is (Figure (\ref{fig:james})) :
\begin{itemize}
 \item Step 1: Solve the Poisson equation on domain $D_1$ with homogeneous Dirichlet boundary conditions
\begin{equation}
 -\Delta \phi_1 = \rho \quad on \quad D_1, \quad \phi_1 = 0 \quad on \quad \partial D_1.
\end{equation}
 \item Step 2: Calculate the surface charge on $\partial D_1$
\begin{equation}
 \partial \rho_1 = \frac{\partial \phi_1}{\partial \boldsymbol{n}} \quad on \quad \partial D_1.
\end{equation}
 \item Step 3: Perform a boundary to boundary convolution from $\partial \rho_1$ to $\partial D_2$
\begin{equation}
 \partial \phi_2 = \int_{\partial D_1} G(\boldsymbol{x}-\boldsymbol{y}) \partial \rho_1(\boldsymbol{y}) dA_{\boldsymbol{y}}
\end{equation}
using fast multiple methods \cite{Greengard87}.
 \item Step 4: Solve another Poisson equation on domain $D_2$ with inhomogeneous Dirichlet boundary conditions
\begin{equation}
 -\Delta \phi_2 = \rho \quad on \quad D_2, \quad \phi_2 =  \partial \phi_2 \quad on \quad \partial D_2.
\end{equation}
\end{itemize}

\subsection{Particle Remapping}
The convergence of particle methods for the one-dimensional Vlasov-Poisson equations has been investigated by Cottet and Raviart \cite{Cottet84}. Their result shows that
particle overlapping and regularization are important for the convergence of the methods. Specifically, the truncation error of a particle method is amplified by a 
time dependent exponential term. Based on Cottet and Raviart's work, we extend the error analysis to PIC methods \cite{Wang10}. Our result is one order higher in the truncation error.
However, as in Cottet and Raviart's analysis, the truncation error is amplified by a time dependent exponential term. The analysis motivates the use of remapping technique, 
a widely used strategy in particle methods in fluid dynamics, to control the exponential error. The basic idea of remapping is simple. 
Since particles will gradually move away from the exact trajectories due to numerical error, we can reduce the displacement by periodically reproducing the distribution function 
$f(\boldsymbol{x}, \boldsymbol{v}, t)$ on a grid by interpolation. A new set of particles, which are created from the grid representation, then replace 
the distorted particle distribution. The later step is identical to the initial step of PIC methods that we initialize the particle positions and weights in equation 
(\ref{eqn:pointparticle}). The error due to remapping will depend on the order of the interpolation function. 
 
In the previous work \cite{Wang10}, we successfully applied the remapped PIC method to the one-dimensional Vlasov-Poisson system. The remapping scheme
was extended in three aspects compared with the standard scheme. First, we used high-order interpolation functions which improve accuracy but
do not preserve positivity. Second, we preserved the positivity of a high-order interpolation by redistributing the 
excess charge into its local neighborhood. The local redistribution algorithm is based on the mass redistribution idea of Chern and Colella \cite{Chern87},
which is first applied to enforce positivity preservation by Hilditch and Colella \cite{Hilditch97}. 
Third, instead of reinitializing on a uniform grid, we reproduced the distribution function on a hierarchy of locally-refined grids. Remapping on a hierarchy of locally-refined 
grids significantly reduces the number of small-strength particles located at the tail of the distribution function. The high-order, positive, and adaptive remapping scheme in 
high dimensional phase space is described below.

\subsubsection{High-order Remapping}
The overall accuracy introduced by remapping will be one order lower than the order of the 
interpolation function since we lose one order of accuracy in the evolution step. For example, the interpolation function with second-order accuracy only results in a first-order
method overall. In this paper, we consider an interpolation function with third-order accuracy derived by Monaghan \cite{Monaghan85}. 
The function in one-dimensional can be expressed as
\begin{equation}
W_4(x, h) = \begin{cases}
1-\frac{5s^2}{2}+\frac{3s^3}{2} & 0 \le s=\frac{|x|}{h} \le 1 \\
\frac{1}{2}{(2-s)}^2(1-s) & 1 \le s=\frac{|x|}{h} \le 2 \\ 
0 & \mbox{otherwise.} \\
\end{cases}
\label{eqn:M4'}
\end{equation}
The one-dimensional expression can be generalized to four dimensions by tensor product,
\begin{equation}
 \boldsymbol{W}_4(\boldsymbol{x}_{\boldsymbol{i}}-\boldsymbol{x}_{k}) = 
\prod_{d=0}^3 W_4(x^d_{\boldsymbol{i}} - x_{k}^d, h^d),
\end{equation}
where $h^d$ is the remapping mesh spacing in phase space. $\boldsymbol{i}$ and $k$ denote the index for the cell-centered grid and the particles, respectively.   

This function, called a modified B-spline, conserves the total charge and represent a quadratic polynomial exactly. In addition, the first- and the second-order derivative of $W_4(x,h)$ are 
continuous. The smoothness property of this modified B-spline is particularly good for scattered data interpolation. However, as with all other high-order interpolation functions, $W_4$ is not positivity preserving.
An interpolation function without positivity might create nonphysical negative charge. This should be avoided in simulations.

\subsubsection{Positivity}
The positivity preserving algorithm is based on the mass redistribution idea
of Chern and Colella \cite{Chern87}, first applied to enforce positivity preservation by Hilditch and Colella \cite{Hilditch97}.
In the algorithm, we redistribute the undershoot of cell $\boldsymbol{i}$
\begin{equation}
\delta f_{\boldsymbol{i}} = \min(0, f_{\boldsymbol{i}}^n) 
\end{equation}
to its neighboring cells $\boldsymbol{i+\ell}$ in proportion to their capacity $\xi$
\begin{equation}
  \xi_{\boldsymbol{i+\ell}} = \max(0, f_{\boldsymbol{i+\ell}}^n) .
\end{equation}
The distribution function is conserved, which fixes the
constant of proportionality
\begin{equation}
 f_{\boldsymbol{i+\ell}}^{n+1} = f_{\boldsymbol{i+\ell}}^n + \frac{\xi_{\boldsymbol{i+\ell}}}{\sum\limits_{\boldsymbol{k}\ne 0}^{\rm neighbors} \xi_{\boldsymbol{i+k}}} \delta f_{\boldsymbol{i}}
\end{equation}
for $\boldsymbol{\ell}\ne 0$ such that cell $\boldsymbol{i+\ell}$ is a neighbor of cell $\boldsymbol{i}$.
Superscript $n$ and $n+1$ denote the interpolated value before and after redistribution, respectively.

The drawback of this approach is that positivity is not guaranteed in a
single pass.  One might have to apply the method iteratively.  In practice, however, we find a few iterations are sufficient.

\subsubsection{Mesh Refinement}
Mesh refinement is an attractive option in improving the efficiency of phase-space remapping. The distribution function in phase space is inhomogeneous,
for example, a Maxwellian distribution in velocity space. When we represent the system by particles, it is best that we can have each particle carries 
a similar amount of weights. Remapping on a hierarchy of locally-refined grids is a good strategy for creating a such set of particles. From another point of view,  
remapping through interpolation is a numerically diffusive procedure. This results in a large number of small-strength particles near the 
tail of the distribution function. The situation becomes worse as we apply remapping frequently. Remapping on a hierarchy of locally-refined
grids, with a coarser grid covering the tail of the distribution function, can reduce the number of those small-strength particles. In the following,
we present the algorithm of remapping with mesh refinement. In designing the algorithm, we have two guiding principles: the total charge 
should be conserved and the overall accuracy on the field needs to be maintained. 

\begin{figure}
\begin{center}
\resizebox{5cm}{!}{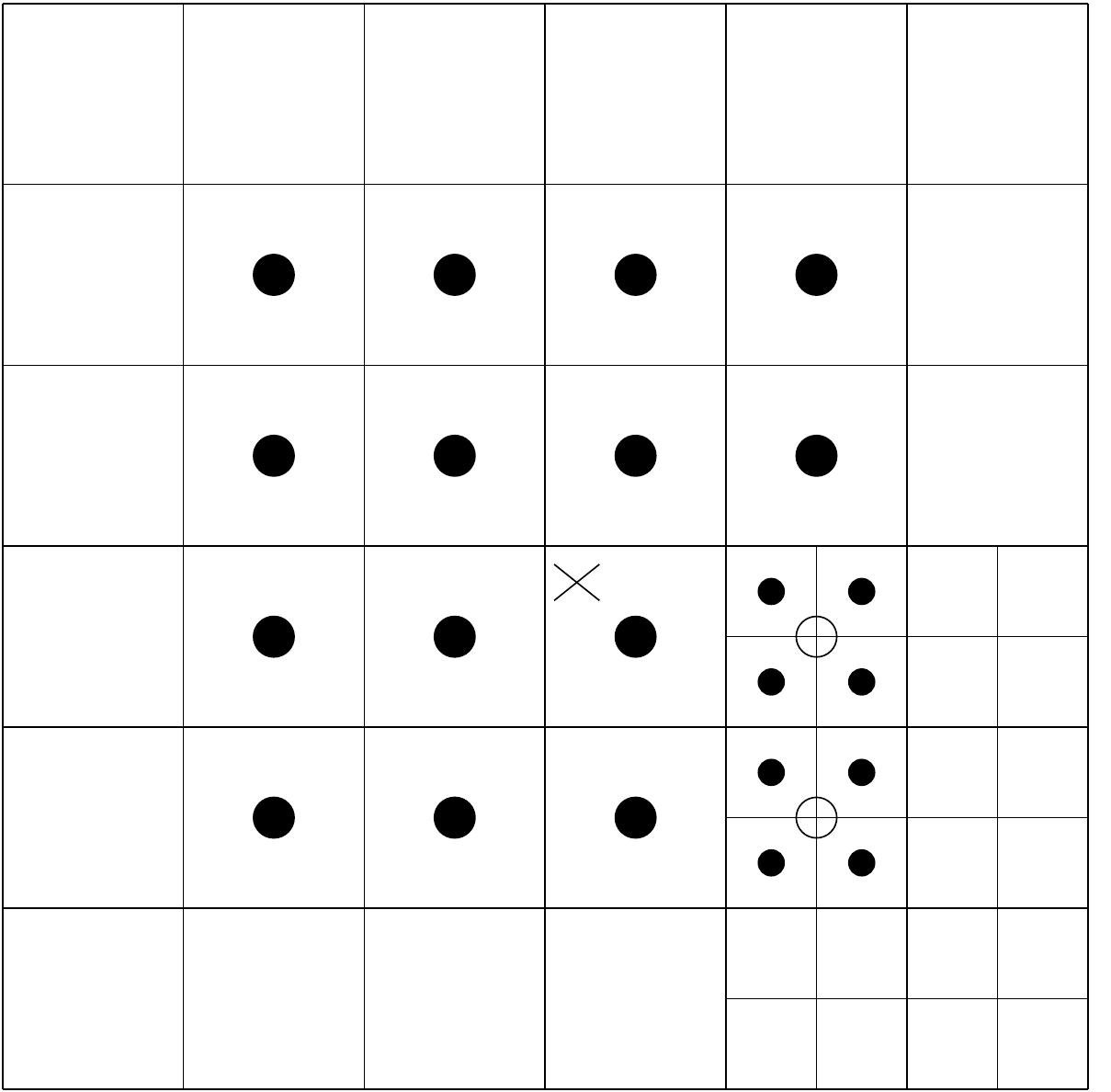} \quad \quad \quad
\resizebox{5cm}{!}{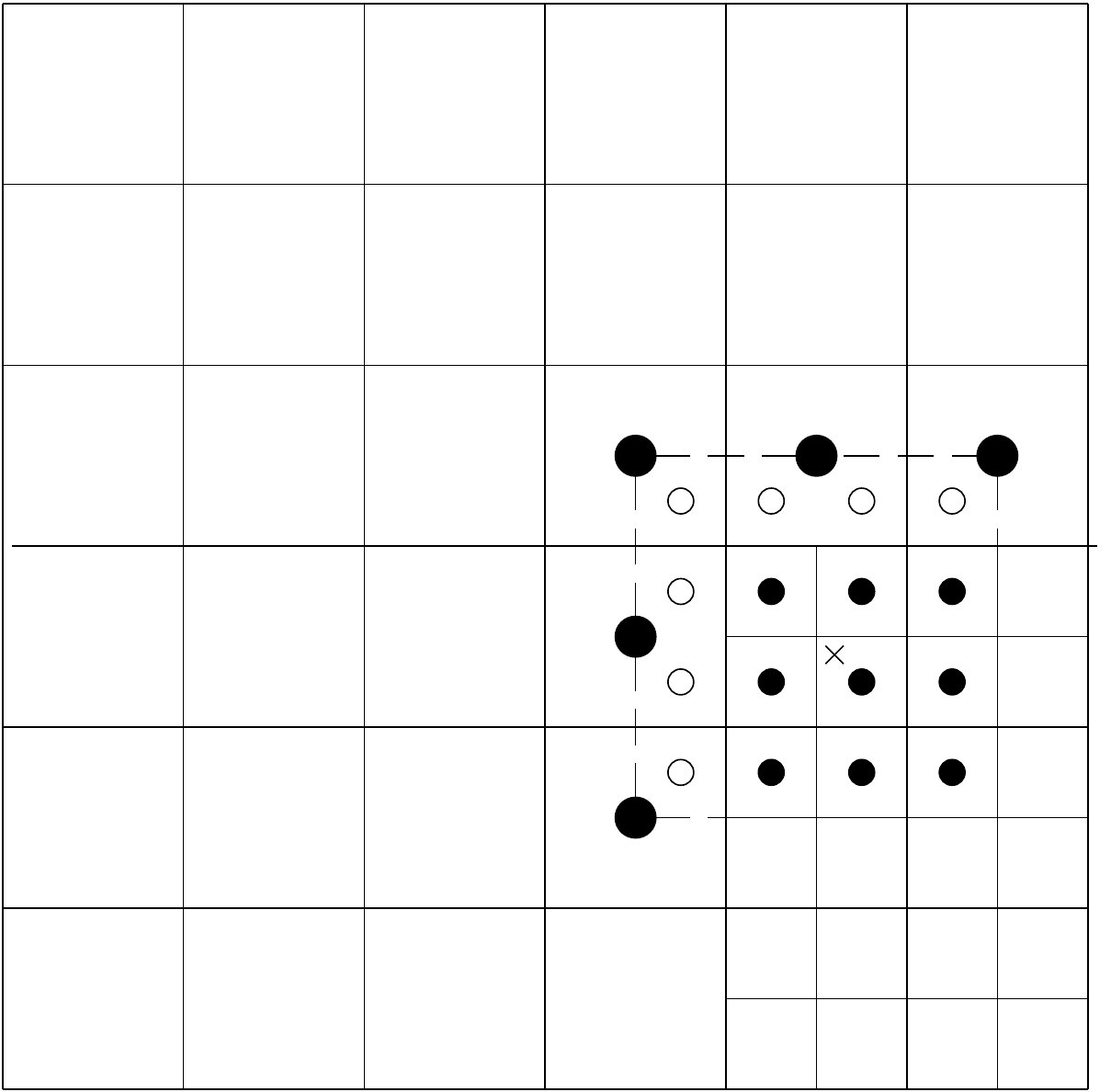}
\caption{ Cross signs denote the particle locations. The valid and the invalid deposited cells are denoted by 
filled circles and open circles, respectively. The refinement ratio is $\boldsymbol{r}_0 = (2, 2)$ in the plots.
Left: Particle is at the coarser level side. Right: Particle is at the finer level side. Cross signs denote the particle locations. \label{fig:mlremap} } 
\end{center}
\end{figure}

Before explaining the algorithm, we introduce the definition of a composite grid. 
We define a hierarchy of cell-centered grids $\Omega_{\ell}$, where $0 \leq \ell \leq \ell_{\rm max}$. $\Omega_{\ell=0}$ is the coarsest grid that covers the whole
 problem domain. The finer grids $\Omega_{\ell>0}$ are constructed as a union of cell-centered rectangles (see Figure \ref{fig:mlremap}). The mesh spacing of each level is 
$\boldsymbol{h}_\ell = \boldsymbol{h}_{\ell-1}/\boldsymbol{r}_{\ell-1}$, where $\boldsymbol{r}_{\ell-1}$ is the refinement ratio of level $\ell-1$.
In four dimensions, $\boldsymbol{h}_\ell \in \mathbb{R}^4$ and  $\boldsymbol{r}_\ell \in \mathbb{Z}^4$.
 The composite grid consists of valid grids
at all levels, where a valid grid is defined as a region not overlain by a finer grid. That is,
\begin{equation}
 \Omega_c = \sum_{\ell=0}^{\ell_{\max}} \left( \Omega_{\ell} \setminus P_{\ell+1}^{\ell}(\Omega_{\ell+1}) \right),
\label{eqn:compgrid}
\end{equation}
where $P_{\ell+1}^{\ell}$ is the operator projecting from level $\ell+1$ to level $\ell$.

At the beginning, a set of particles are created from the cell center of the composite grid. In the remapping step, each particle first finds the valid cell 
in the composite grid it belongs to. One particle can only belong to a single valid cell.
If the cell is far enough away from a coarse fine-interface, the charge can be interpolated on the grid as in equation (\ref{eqn:M4'}).
If the cell is near a coarse-fine interface such that the interpolation stencil 
intersects the coarse-fine interface, special care must be taken. First, we interpolate the charge on the surrounding cells as usual. After
deposition, we know that not all deposited cells are valid cell. We need a further step to transfer the charge from invalid cell to valid cell. There are two 
cases depending on where the invalid cell is located. If the invalid cell is in a coarser level and it is covered by the valid cells of a finer level, we transfer the deposited 
charge from the coarser level to the finer level through interpolation. On another hand, if the invalid cell is outside the grid of the current level, 
the charge is transfered by projection. Figure \ref{fig:mlremap} shows the algorithm in two dimensions. The four-dimensional case can be generalized easily.
It is worth mentioning that we lose one order of accuracy in interpolating the coarser level charge into the finer level. However, since the coarse-fine 
interface is in co-dimension one, the expected accuracy in the field, e.g., second-order, will be preserved in $L_\infty$ norm error. 
Our current implementation doesn't have time-dependent adaptivity. This feature can be incorporated by selecting some refinement criterion, for example each cell in phase space has
similar number of particles.

\section{Parallel Implementation and Issues}
\label{sec:implementation}
The parallel implementation of the algorithm is straightforward based on domain decomposition in physical space. The physical space is decomposed 
into $M$ patches. Given a parallel machine with $N$ processors, each patch is assigned to a processor cyclically. 
Particles are assigned to patches according to their physical space positions.
Using MPI, patches communicate with each other through ghost-cells and particles move between patches. 
This is the default implementation in Chombo software \cite{chombo}.

The current implementation is not a scalable algorithm (weak scaling) because of decomposition in physical space only.
The potential issue is that 
when the problem size increases, since the number of processors is scaled in proportion to the problem size in physical space, the computational time and memory usage
will increase in proportion to the problem size in velocity space. In the worst case, the processor will be out of memory. An alternative implementation is 
based on domain decomposition in phase space. In this implementation, it might be the case that  
particles belong to the same cell in physical space (equation (\ref{eqn:rho})) are distributed on different processors. Since the Laplacian operator is linear, we can choose to solve 
the Poisson equation separately on different processors. The total fields are then obtained by MPI\_Allreduce.


\section{Numerical Tests}
\label{sec:tests}
We demonstrate PIC methods with adaptive phase-space remapping on a set of classical plasma and beam problems in two dimensions,
including linear Landau damping, the two stream instability, and beam propagation in the paraxial model. 
For satisfying the overlapping condition, we choose $\varepsilon_x/h_x = \varepsilon_y/h_y = 2$, where $\boldsymbol{\varepsilon_x} = (\varepsilon_x, \varepsilon_y)$.

We use Richardson extrapolation for error estimate. If $\tilde{\boldsymbol{E}}^h$ is the electric field computed with the initial phase space discretization 
$(h_x, h_y, h_{v_x}, h_{v_y})$ and integration step size $\Delta t$, and $\tilde{\boldsymbol{E}}^{2h}$ computed with $(2h_x, 2h_v, 2h_{v_x}, 2h_{v_y})$ and 
$2\Delta t$, the relative solution error in direction $d$ is defined as
\begin{equation}
 e^{h}_d = |\tilde{E}^h_d - \tilde{E}^{2h}_d|.
\end{equation}
 $q$ is the order of the method and is calculated by
\begin{equation}
 q= \min_d \log_2 \left( \frac{||e^{2h}_d||}{||e^{h}_d||}\right).
\end{equation} 

\subsection{Linear Landau Damping}
The initial distribution for linear Landau damping is 
\begin{equation}
f_0(x, y, v_x, v_y)  = \frac{1}{2\pi} \exp(-(v_x^2+v_y^2)/2)(1+\alpha \cos(k_xx)\cos(k_yy)), \\
\label{eqn:linearlandau}
\end{equation}
where $\alpha=0.05$, $k_x=k_y=0.5$ and $v_{\rm max}=6.0$. The physical domain is $(x, y)\in [0, L=2\pi/k_x] \times [0, L=2\pi/k_y]$ with periodic boundary 
conditions. Particle charges with strength less than $1.0\times10^{-9}$ are ignored. In the simulation, we apply remapping every 5 PIC time steps.

In the first test, we are interested in the evolution of the amplitude of the electric field.
According to Landau's theory, the electric field is expected to decrease exponentially with damping rate $\gamma=-0.394$. The behavior of exponential decay has been 
observed by many other authors, mostly calculated by grid methods \cite{Nakamura99,Filbet01,Crouseilles08II,Crouseilles09}.

\begin{figure}[htp]
\centering
  \subfigure[with remapping]{\includegraphics[scale=0.8]{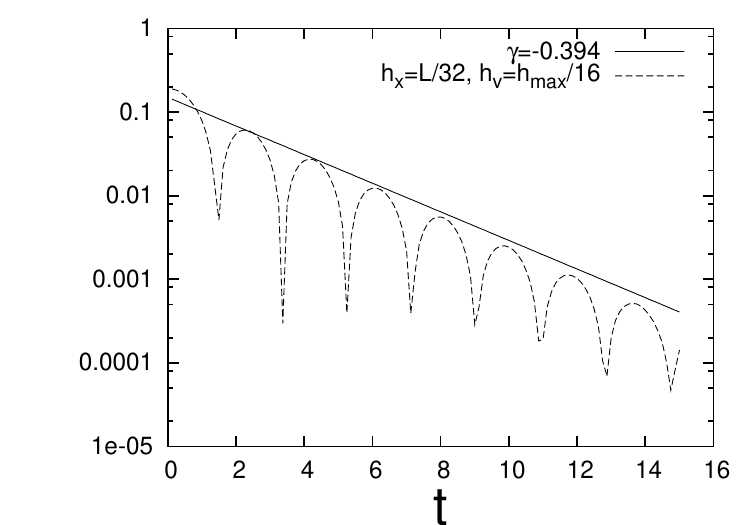}} 
  \subfigure[without remapping]{\includegraphics[scale=0.8]{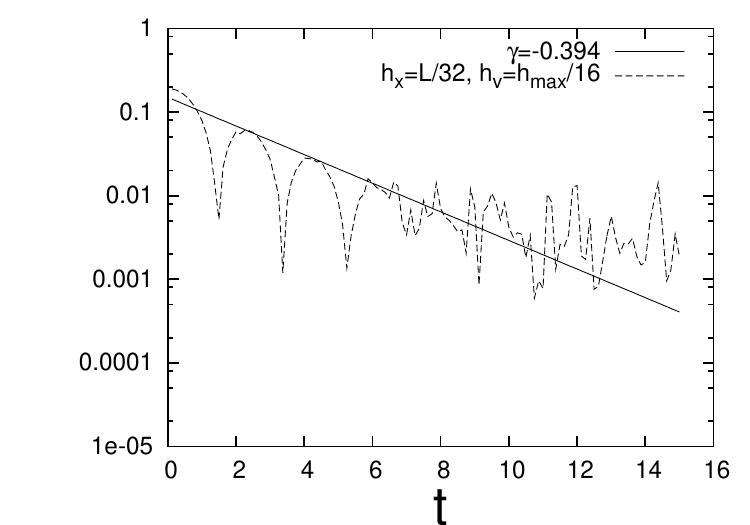}}
\caption{ The amplitude of the electric field for 2D linear Landau damping problem. 
Scales $(h_x, h_{v_x})$ above denote the particle grid mesh spacing at the base level, where $h_x=h_y$ and $h_{v_x}=h_{v_y}$.
With remapping, the computed damping rate is very close to the theoretical value $\gamma=-0.394$.
Without remapping, the simulation fails to track the exponential decay after a few damping circles.
\label{fig:linearlandauenergy} }
\end{figure}

We initialize the problem on two levels of grids with base level at $h_x=h_y=L/32,
h_{v_x}=h_{v_y}=v_{\rm max}/16$. The velocity space is refined on sub-domain $\boldsymbol{v} \in [-3, 3] \times [-3, 3]$
with a refinement ratio $2$. The PIC step size is $dt=1/8$.
We compare the simulation with and without remapping in Figure (\ref{fig:linearlandauenergy}). 
In the case with remapping, the computed damping rate is very close to the theoretical value. 
The simulation without remapping fails to track the exponential decay.

\begin{figure}[htp]
\centering
  \subfigure[$L_\infty$ errors of $\tilde{E}_x$]{\includegraphics[scale=1.0]{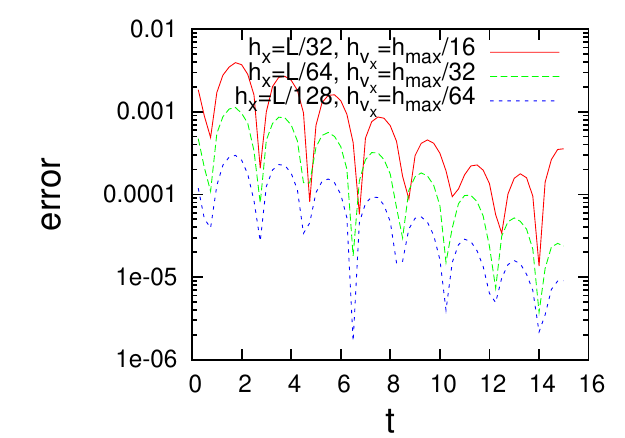}}
  \subfigure[Convergence rates]{\includegraphics[scale=1.0]{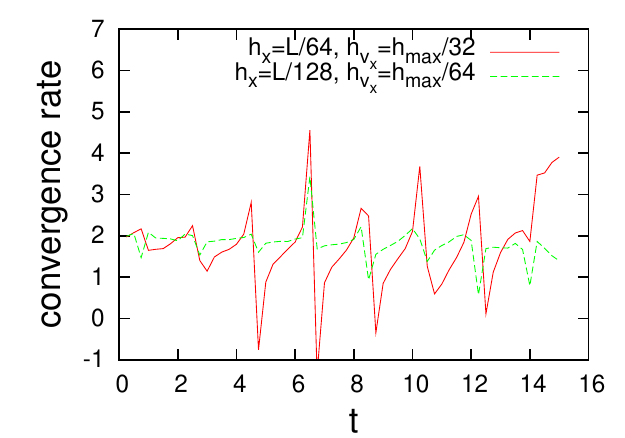}}
\caption{ Error and convergence rate plots for 2D linear Landau damping problem with remapping.
Scales $(h_x, h_{v_x})$ above denote the particle grid mesh spacing at the base level, where $h_x=h_y$ and $h_{v_x}=h_{v_y}$.
Second-order convergence rates are obtained. 
(a) the $L_\infty$ norm of the electric field errors on three different resolutions. 
(b) the convergence rates for the errors on plot (a).
\label{fig:linearlandauerrconwp} }
\end{figure}

\begin{figure}[htp]
\centering
  \subfigure[$L_\infty$ errors of $\tilde{E}_x$]{\includegraphics[scale=1.0]{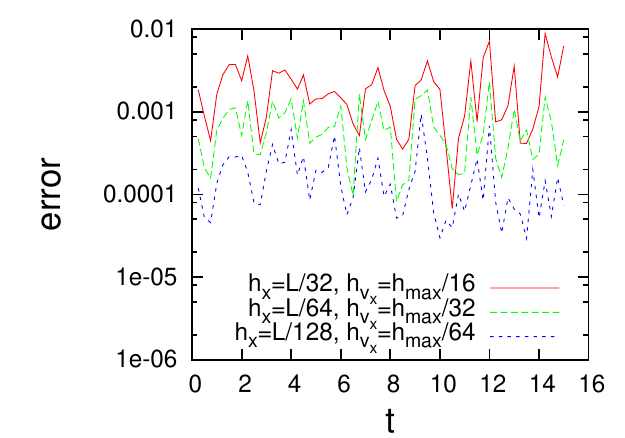}}
  \subfigure[Convergence rates]{\includegraphics[scale=1.0]{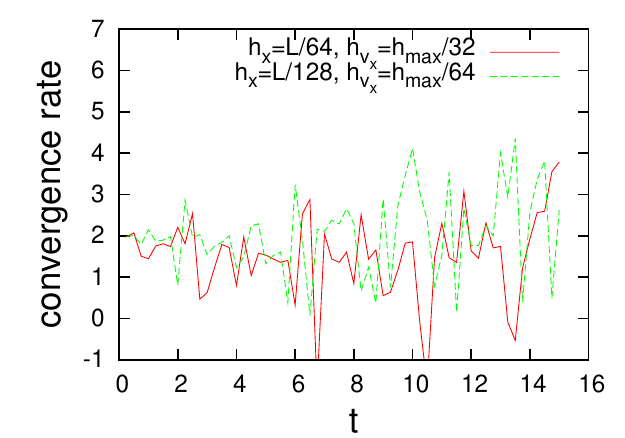}}
\caption{ Error and convergence rate plots for 2D linear Landau damping problem without remapping.
Scales $(h_x, h_v)$ above denote the particle grid mesh spacing at the base level, where $h_x=h_y$ and $h_{v_x}=h_{v_y}$.
The errors without remapping are much larger than the case with remapping (see Figure (\ref{fig:linearlandauerrconwp})).
(a) the $L_\infty$ norm of the electric field errors on three different resolutions. 
(b) the convergence rate for the errors on plot (a).
\label{fig:linearlandauerrconnr} }
\end{figure}

In the second test, we compare the electric field errors and
corresponding convergence rates with and without remapping in Figure
(\ref{fig:linearlandauerrconwp}) and (\ref{fig:linearlandauerrconnr}).
We see that remapping significantly reduces the electric field errors and
improves their corresponding convergence rates.

\subsection{The Two Stream Instability}
The initial distribution for the two stream instability is 
\begin{equation}
f_0(x, y, v_x, v_y)  = \frac{1}{12\pi} \exp(-(v_x^2+v_y^2)/2)(1+\alpha \cos(k_xx)) (1 + 5v_x^2), \\
\label{eqn:twostream}
\end{equation}
where $\alpha=0.05$, $k_x=0.5$ and $v_{\rm max}=9.0$. The physical domain is $(x, y)\in [0, L=2\pi/k_x] \times [0, L=2\pi/k_y]$ with periodic boundary 
conditions. Particle charges with strength less than $1.0\times10^{-9}$ are ignored. As in linear Landau damping problem, we apply remapping every
5 PIC time steps. The velocity space is refined on sub-domain $\boldsymbol{v} \in [-4.5, 4.5]\times[-4.5, 4.5]$ with a refinement ratio $2$. 

\begin{figure}[htp]
\centering
  \subfigure[$L_\infty$ errors of $\tilde{E}_x$]{\includegraphics[scale=1.0]{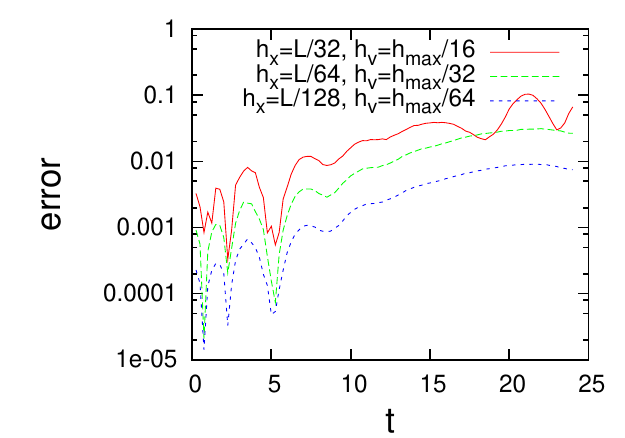}}
  \subfigure[Convergence rates of the errors on the left]{\includegraphics[scale=1.0]{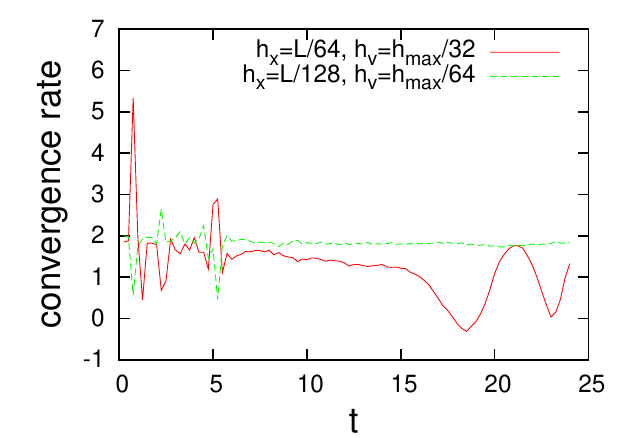}} 
\caption{ Error and convergence rate plots for the two stream instability with remapping.
Scales $(h_x, h_v)$ above denote the particle grid mesh spacing at the base level, where $h_x=h_y$ and $h_{v_x}=h_{v_y}$. The PIC step size is $dt=1/8$ at the lowest resolution.
(a) the $L_\infty$ norm of the electric field errors on three different resolutions. 
(b) the convergence rate for the errors on plot (a).
\label{fig:twostreamerrconwp} }
\end{figure}

\begin{figure}[htp]
\centering
  \subfigure[$L_\infty$ errors of $\tilde{E}_x$]{\includegraphics[scale=1.0]{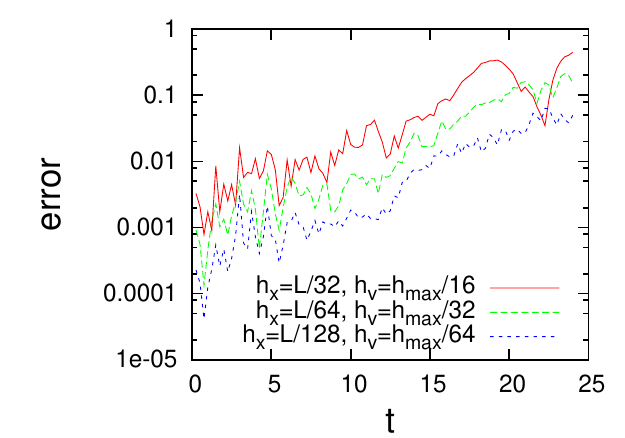}}
  \subfigure[Convergence rates of the errors on the left]{\includegraphics[scale=1.0]{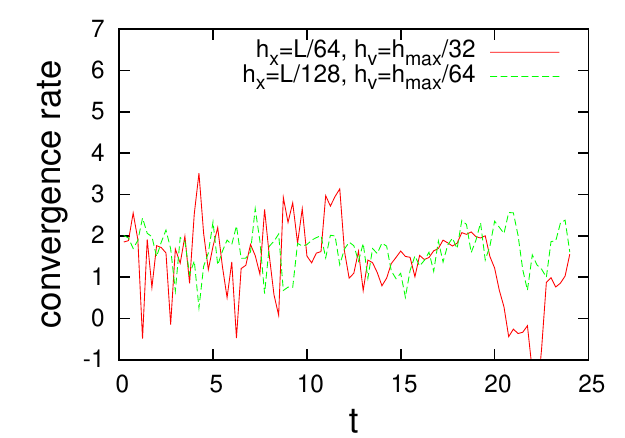}} 
\caption{ Error and convergence rate plots for the two stream instability without remapping.
Scales $(h_x, h_v)$ above denote the particle grid mesh spacing at the base level, where $h_x=h_y$ and $h_{v_x}=h_{v_y}$. The PIC step size is $dt=1/8$ at the lowest resolution.
(a) the $L_\infty$ norm of the electric field errors on three different resolutions. 
(b) the convergence rate for the errors on plot (a). \label{fig:twostreamerrconnr} }
\end{figure}

We compare the electric field errors and their convergence rates
 with and without remapping as before. Figure (\ref{fig:twostreamerrconnr}) shows the $L_\infty$ norm of the errors at the case without remapping in three
different resolutions. The corresponding convergence rates are shown on the right of the error plots. Second-order convergence rates are lost at the 
early time of the simulation. Comparing with the results with remapping in Figure (\ref{fig:twostreamerrconwp}), we see that 
remapping extends the second-order convergence rates to longer times.

\begin{figure}[htp]
\centering
  \subfigure[with remapping at time $t=20$]{\label{fig:twostreamwp} \includegraphics[scale=0.3]{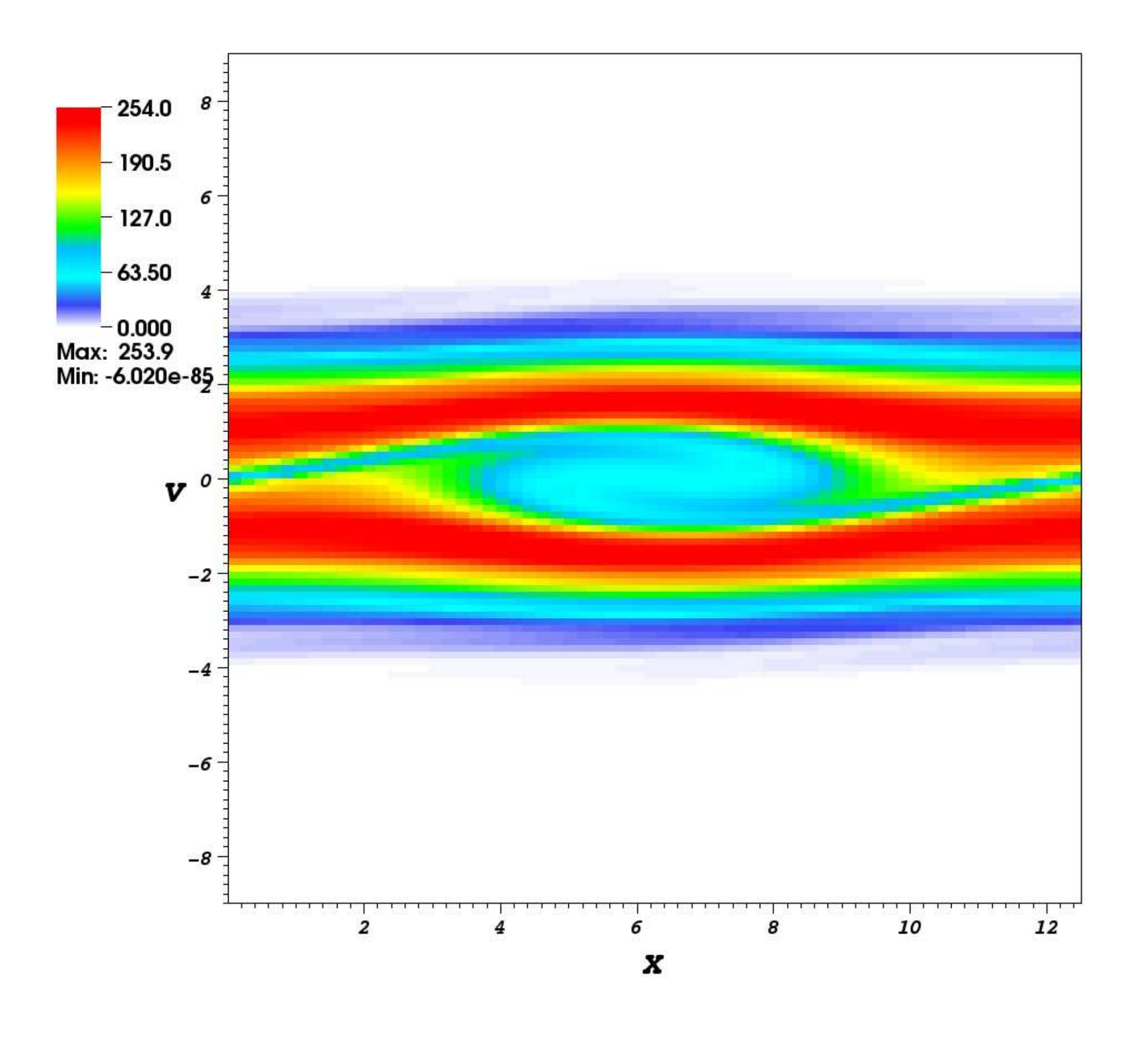}}
  \subfigure[without remapping at time $t=20$]{\label{fig:twostreamnr} \includegraphics[scale=0.3]{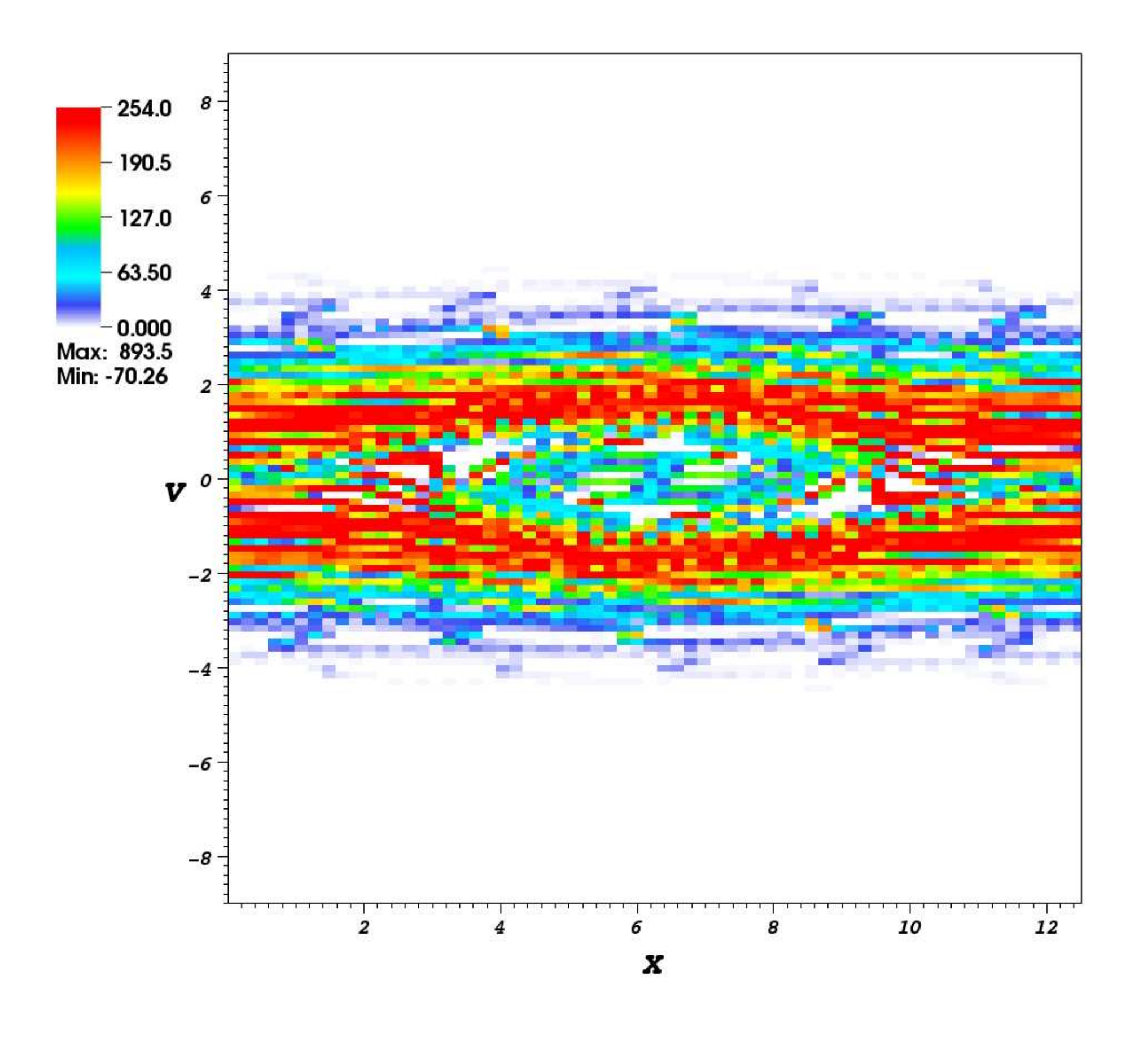}} 
\caption{ Comparison of $F(x,v_x)$, the distribution function projected on space $(x, v_x)$, at the same instant
  of time $t=20$ with (Left) and without remapping (Right) for the two stream instability problem. The projected value is 
$F(x, v_x) = \int_{0}^{L} \int_{-\infty}^{\infty} f(x, y, v_{x}, v_{y}) dy dv_y$. The grid-based distribution function is obtained by reproducing the 
particle-based distribution function through a second-order interpolation. We initialize the distribution function on two levels of grids, with base
level at $h_x = h_y = L/64$, $h_{v_x} = h_{v_y} = v_{\rm max}/32$. The grid is refined by factor of 2 in velocity space on sub-domain 
$\boldsymbol{v} \in [-4.5, 4.5]\times[-4.5, 4.5]$. The classical PIC method results in a noisy solution with large errors in maximum. 
Both numerical noise and errors in maximum are significantly reduced by using remapping. 
The negative minimum in the case without remapping is a superficial effect due to project 4D data to 2D using a high-order interpolation for visualization purpose.
\label{fig:twostreamcompdist} }
\end{figure}

We also compare the projected distribution function on plane $(x, v_x)$ at the same instant time $t=20$ by both methods in Figure 
(\ref{fig:twostreamcompdist}). For visualization purpose, in the case without remapping, we interpolate the particle-based distribution function
to a grid in phase space. We see that the classical PIC method results in a noisy solution (see Figure (\ref{fig:twostreamnr})). 
Figure (\ref{fig:twostreamwp}) shows the distribution function computed by the PIC method with remapping. 

\subsection{Semi-Gaussian Beam}
The paraxial model is an approximation to the steady-state Vlasov-Maxwell equation in three dimensions. The K-V distribution is a measure solution of the paraxial model. 
Given an arbitrary initial distribution, we can focus a beam with the same matching forces for the K-V beam using the concept of equivalent beam.  
Here, we consider an initial semi-Gaussian beam focused by an uniform electric field using the concept of equivalent beam.
The model has been considered by many authors \cite{Sonnendrucker04,Crouseilles08II,Crouseilles09}.

In the test, the beam is composed of ionized potassium. The physical parameters are the following: current $I=0.2A$, beam velocity $v_b = 0.63 \times 10^6 m/s$,
and the radius of the beam $a=0.02m$. We choose the tune depression $\eta = 1/2$. For the normalization of the paraxial model, we refer to the work of Filbet and Sonnendrucker \cite{Filbet06}.
We use normalization parameters $(x_0, v_0) = (a, \frac{\epsilon_x v_b}{2a})$. 
This results in the normalized Poisson system 
 \begin{equation}
 -\bigtriangleup \phi = \int_{\mathbb{R}^2} f d\boldsymbol{v}, \quad -\nabla \phi = \boldsymbol{E},
\label{eqn:beampoisson}
\end{equation}
with initial semi-Gaussian distribution
\begin{equation}
f_{0} (x, y, v_x, v_y) = \left\{
\begin{array}{cl}
\frac{4(1-\eta^2)}{\pi \eta^2}\exp(-(v_x^2+v_y^2)/2), & \mbox{$x^2+y^2 \leq 1$} \\
0, & \mbox{otherwise}
\end{array}
\right.
\end{equation}
and the external matching field
\begin{equation}
 E^{e}(x, y, t) = - \frac{4}{\eta^2} (x\boldsymbol{e}_x + y\boldsymbol{e}_y).
\end{equation}
For the numerical parameters, we choose $(h_x=|L_x|/128, h_y = |L_y|/128)$ and $(h_{v_x}=v_{\rm max}/128, h_{v_y}=v_{\rm max}/128)$, 
where $(L_x, L_y)=(-10, 10)$ and $v_{\max} = 10$. The PIC time step is $dt = 0.00052925$.
 
\begin{figure}[htp] 
  \begin{minipage}[b]{0.45\textwidth}
   \centering 
\includegraphics[width=50mm, height=40mm]{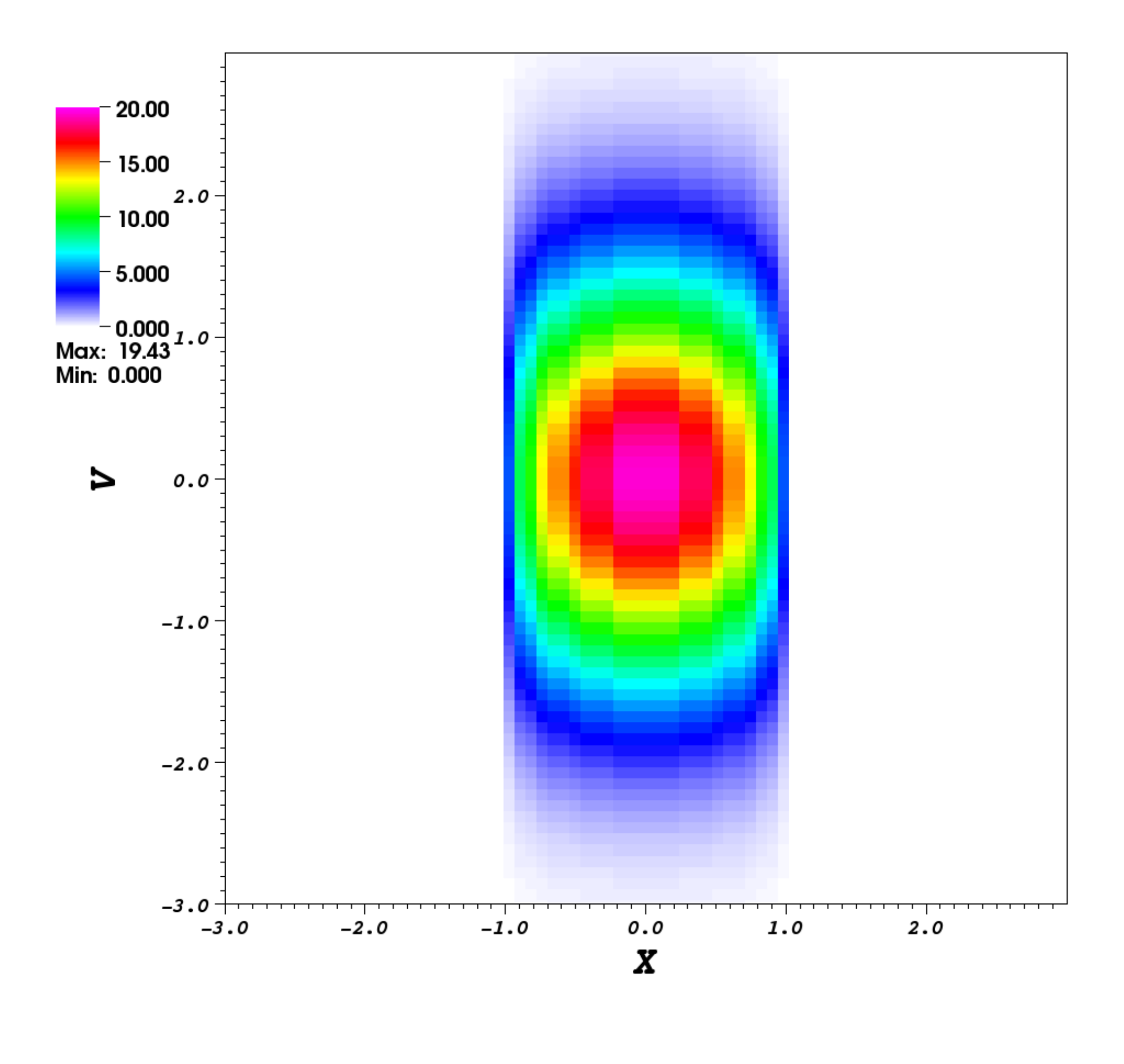}
\includegraphics[width=50mm, height=40mm]{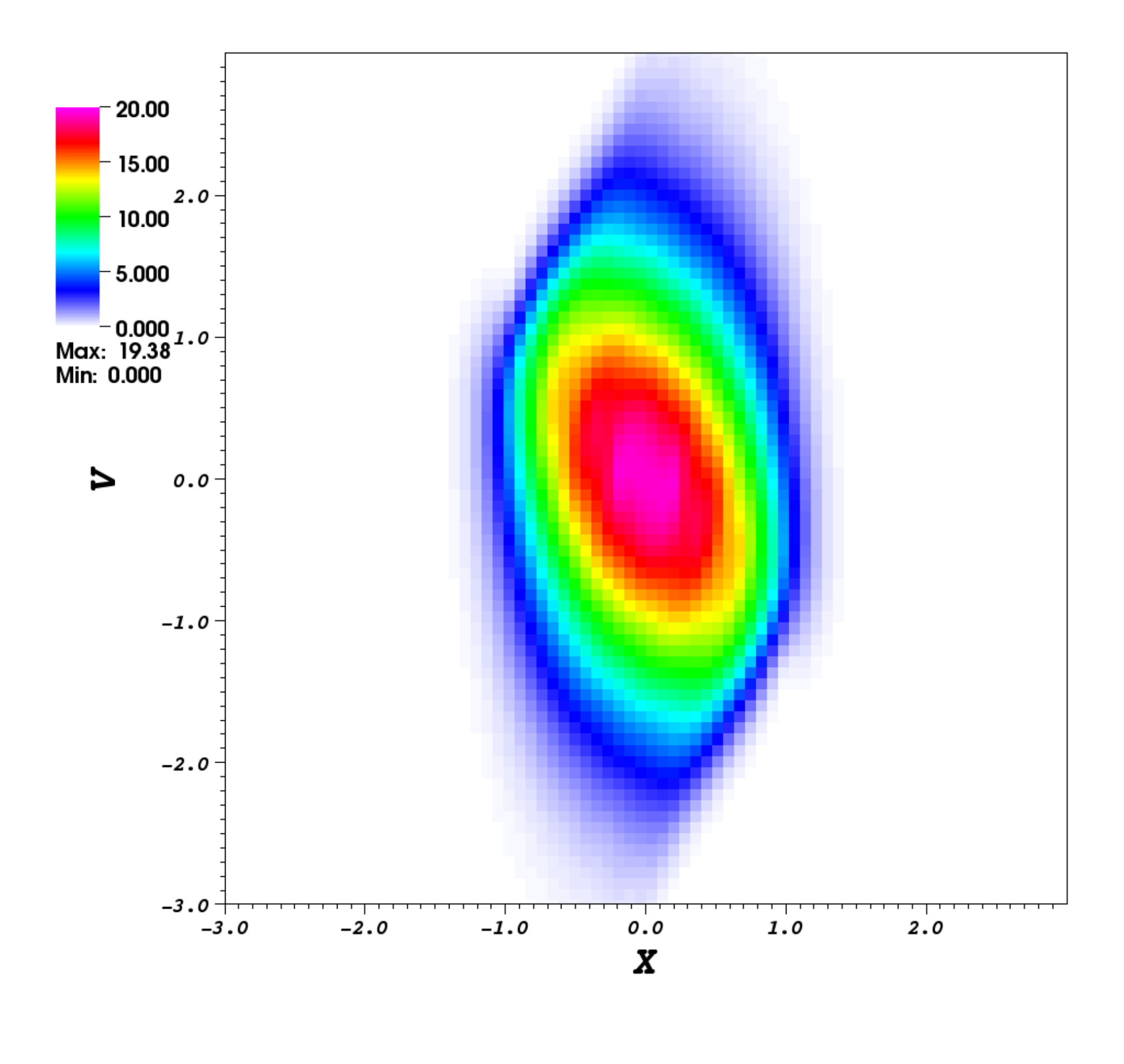}
\includegraphics[width=45mm, height=40mm]{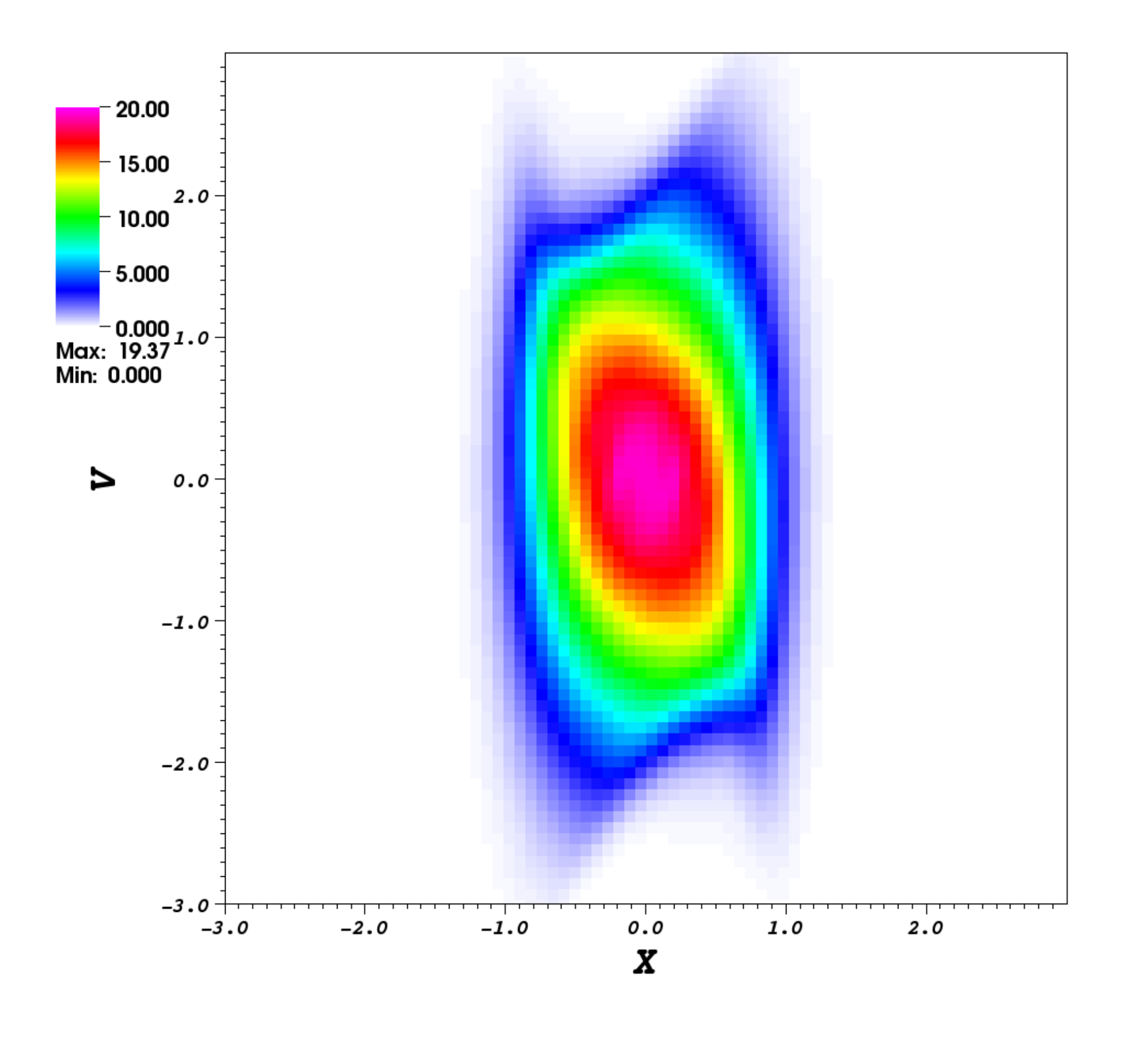}
\includegraphics[width=50mm, height=40mm]{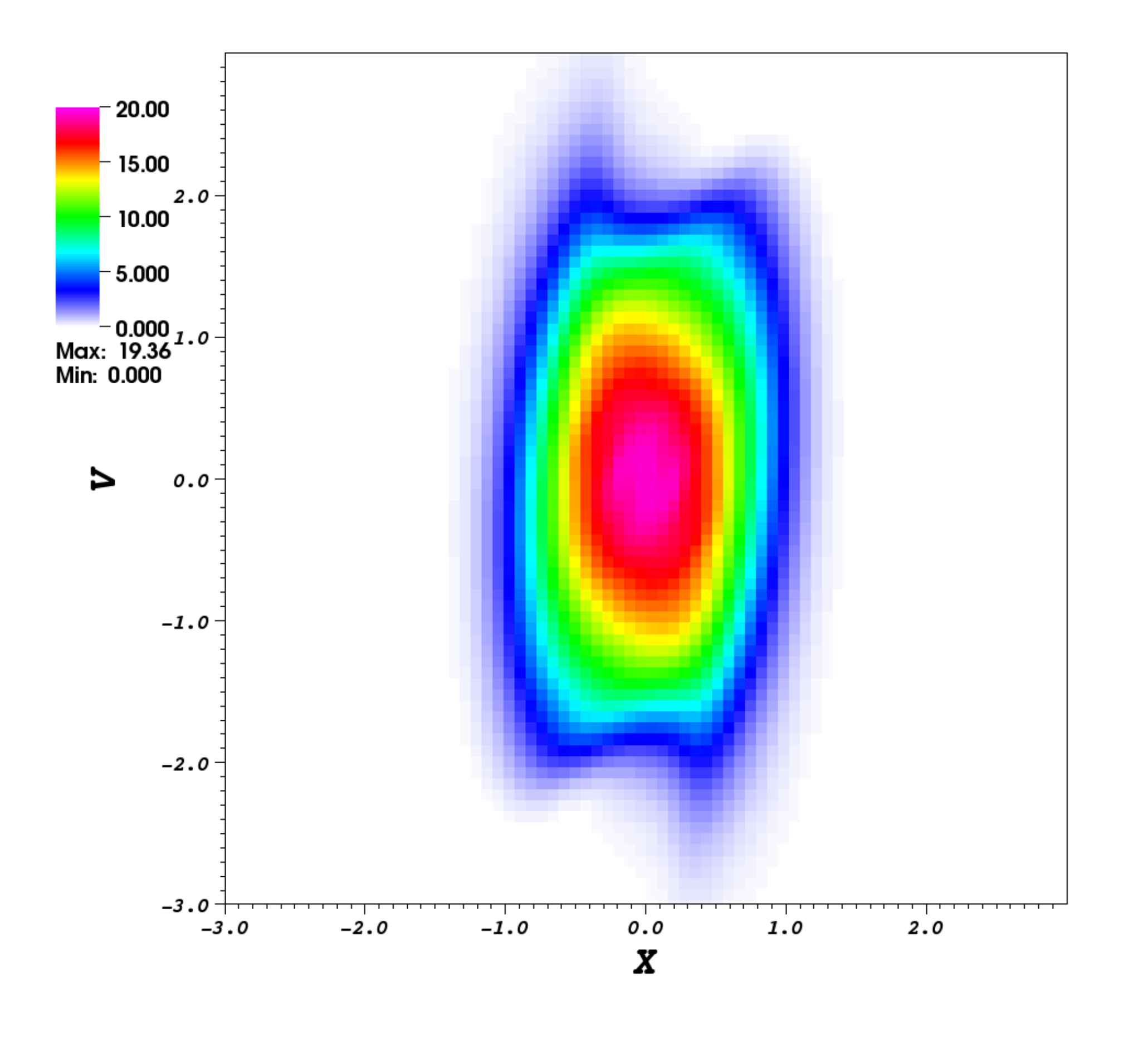}
\includegraphics[width=50mm, height=40mm]{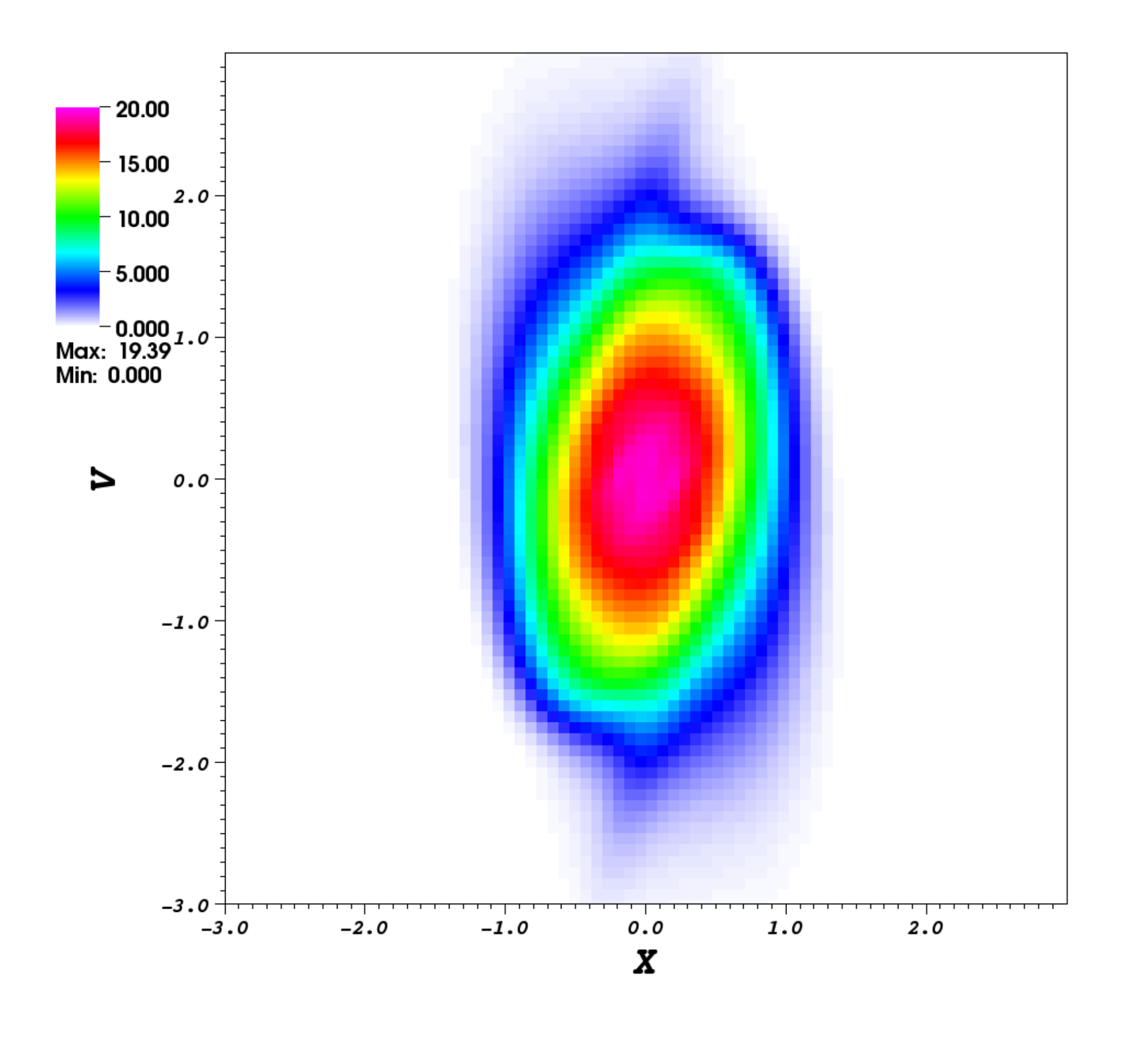}
  \end{minipage}
\hfill
  \begin{minipage}[b]{0.45\textwidth}
   \centering 
\includegraphics[width=50mm, height=40mm]{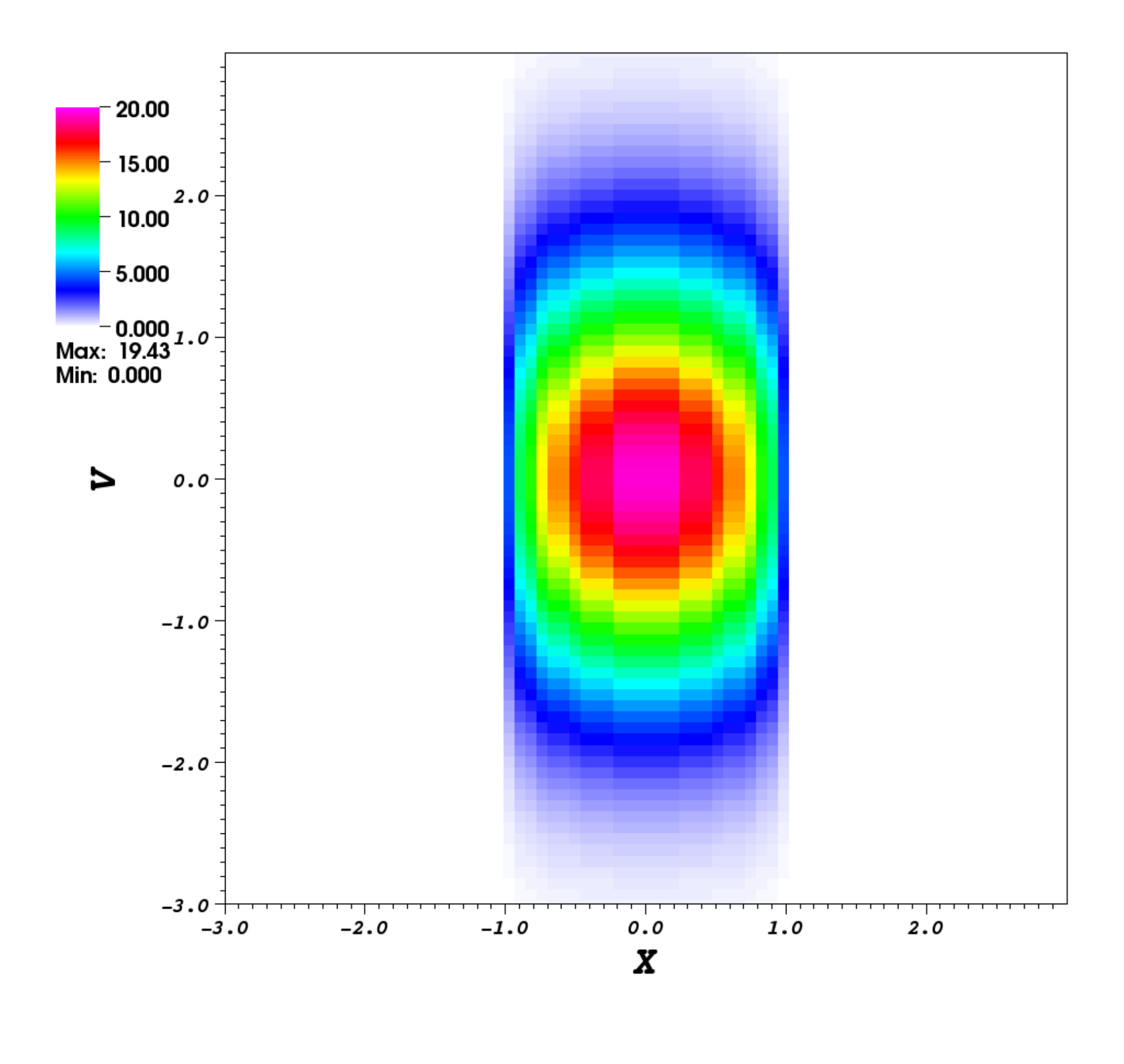}
\includegraphics[width=50mm, height=40mm]{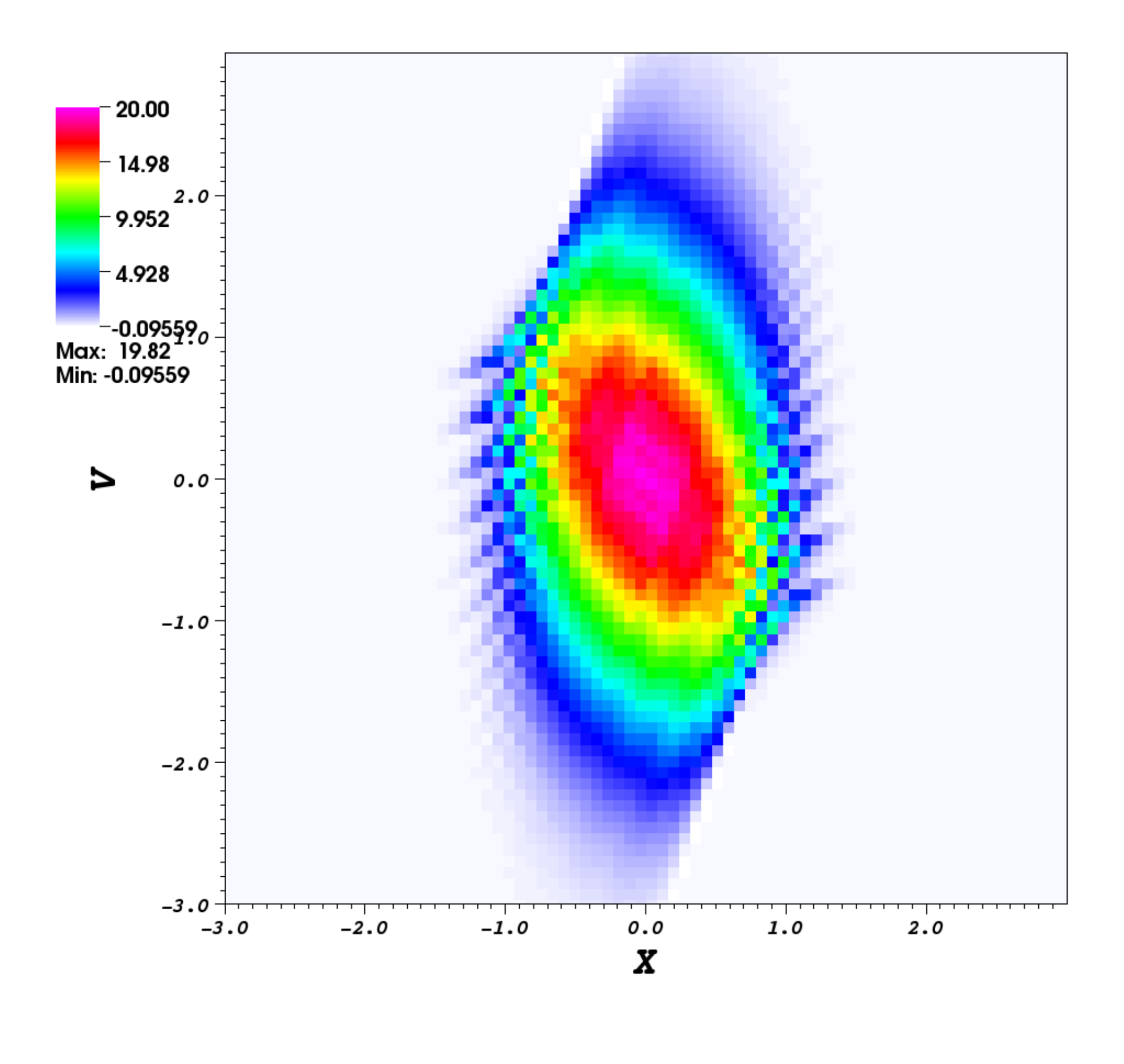}
\includegraphics[width=50mm, height=40mm]{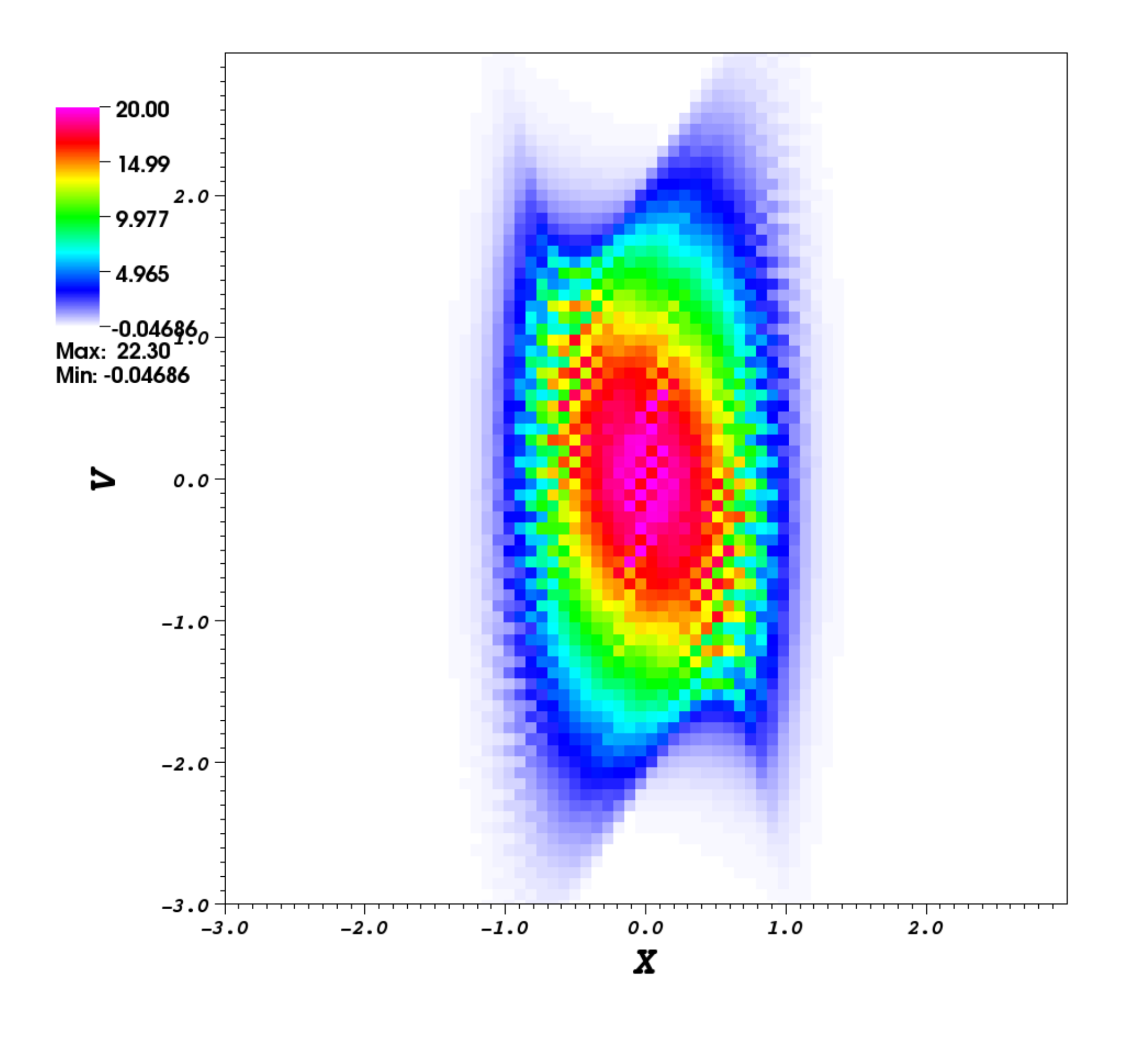}
\includegraphics[width=50mm, height=40mm]{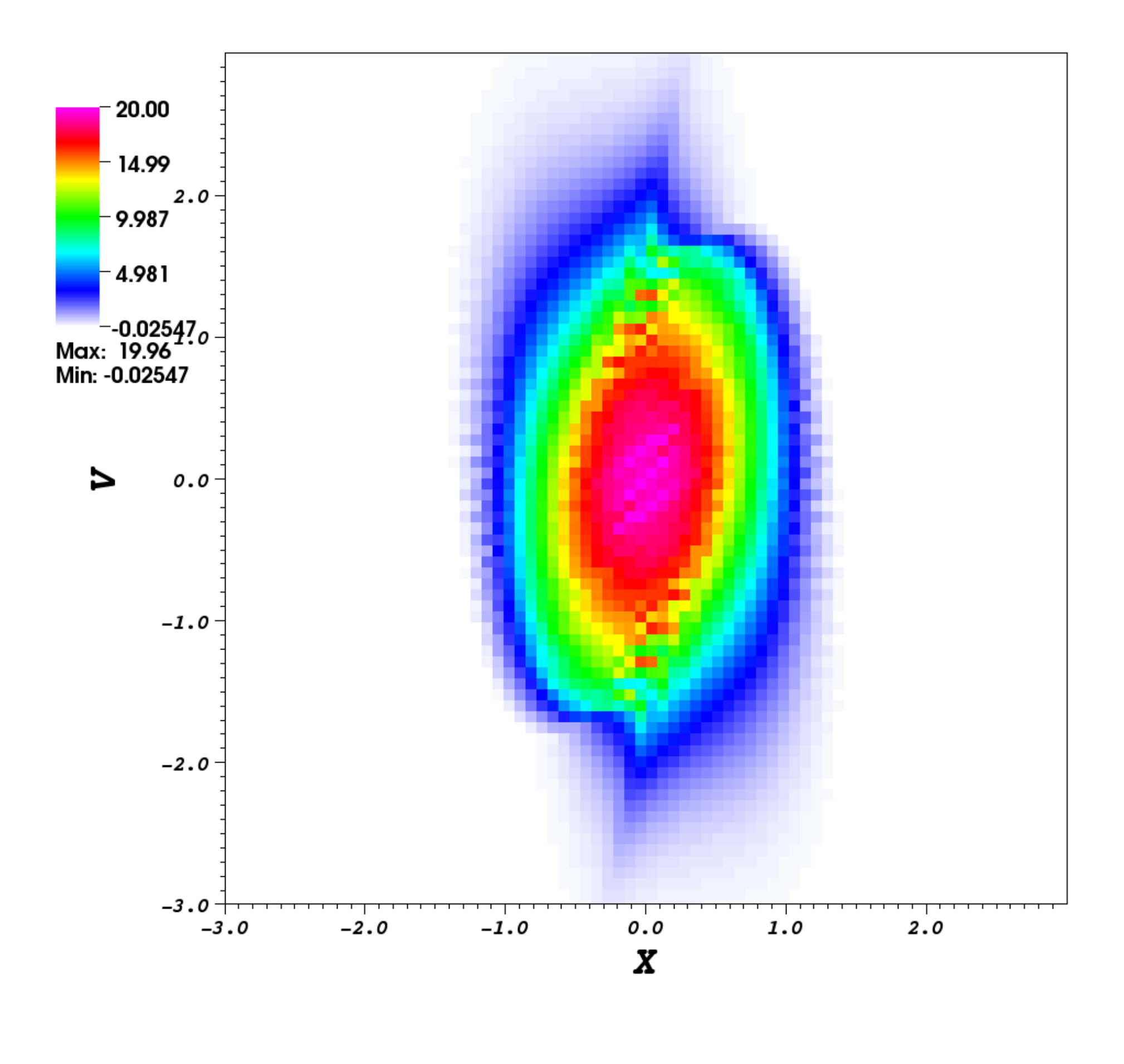}
\includegraphics[width=50mm, height=40mm]{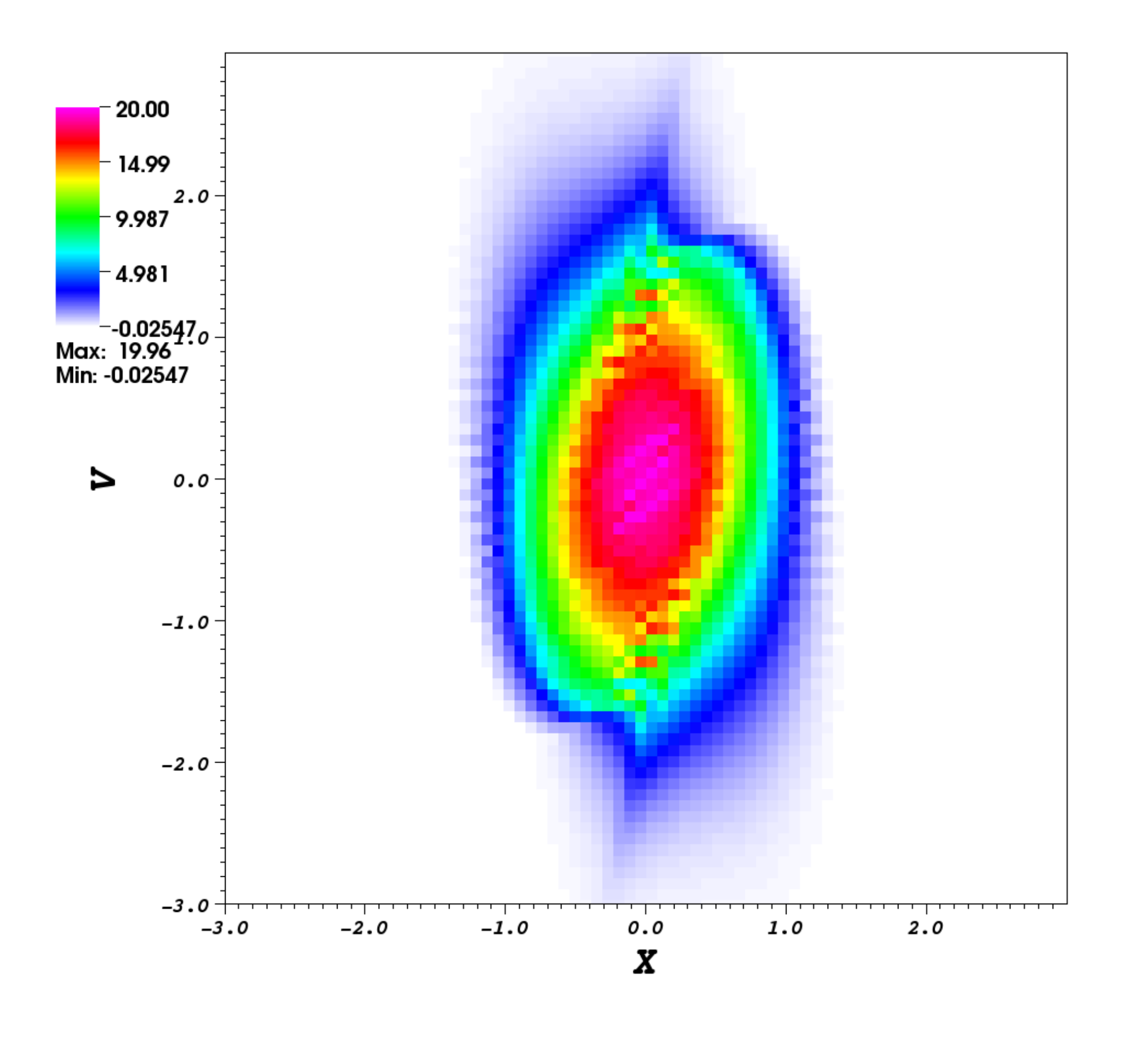}
  \end{minipage}
\caption{Time evolution of the projection of the distribution function on $(x, v_x)$
at time $t=0, 0.3176, 0.5239, 0.7410, 0.9527$ (from top to bottom). 
The projected value is $F(x, v_x) = \int_{0}^{L} \int_{-\infty}^{\infty} f(x, y, v_{x}, v_{y}) dy dv_y$. 
The grid-based distribution function is obtained by reproducing the 
particle-based distribution function through a second-order interpolation.
The columns on the left and on the right are simulations running with and without remapping, respectively.
The classical PIC method results in a noisy solution with large errors in maximum.
\label{fig:2ddistcase8xv}}
\end{figure}
\begin{figure}[htp]
\centering
  \subfigure[RMS of x (64)]{\includegraphics[scale=0.7]{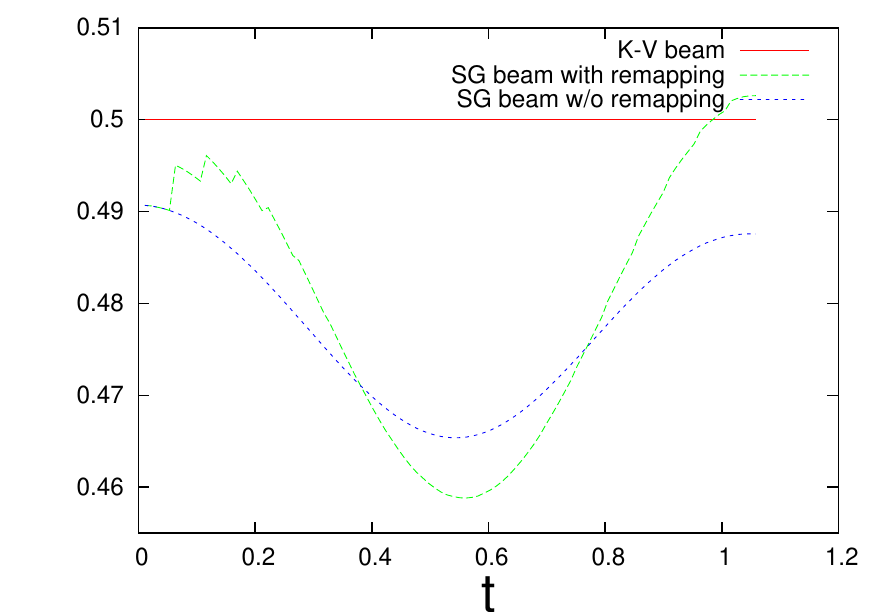}}
  \subfigure[RMS of x (128)]{\includegraphics[scale=0.7]{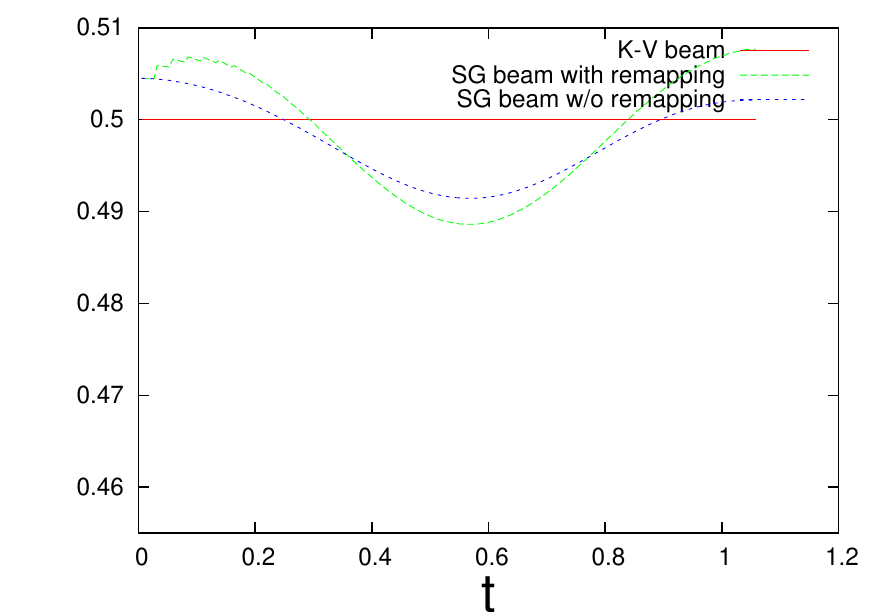}} 
  \subfigure[RMS of $v_x$ (64)]{\includegraphics[scale=0.7]{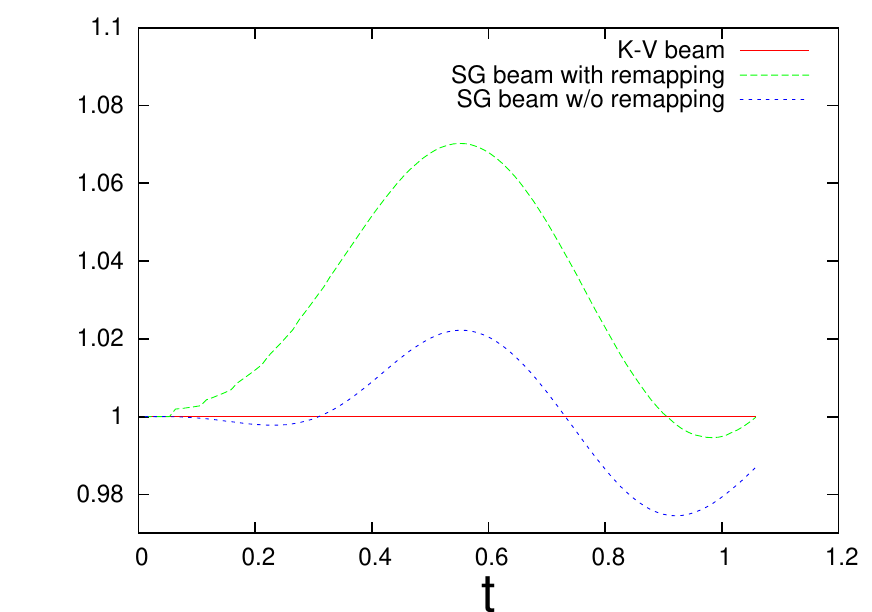}} 
  \subfigure[RMS of $v_x$ (128)]{\includegraphics[scale=0.7]{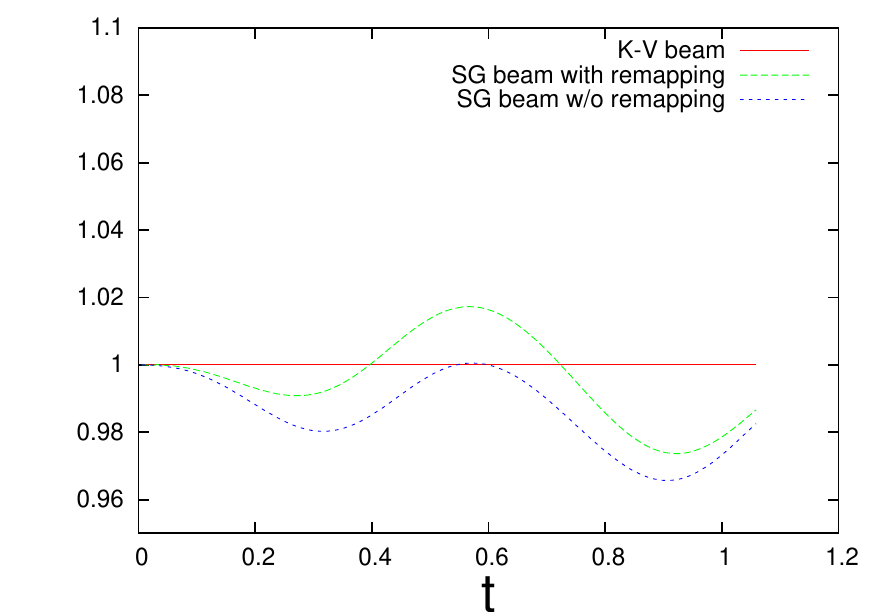}}
\caption{ Comparison of RMS quantities of semi-Gaussian beam with remapping, without remapping, and its 
equivalent K-V beam. We show the quantities in two different resolutions (64 vs 128 grid size in each dimension).
Although the values are oscillatory, they remains close to the equivalent K-V beam and they converge to the equivalent K-V beam as the 
resolution increases. 
\label{fig:rms} }
\end{figure}

Figure (\ref{fig:2ddistcase8xv}) shows the projection of the distribution function on planes ($x, v_x$) with and without remapping, respectively.
The simulation with remapping gives a well-resolved result which preserves the positivity of the distribution function. 
Meanwhile, we show the root mean square (RMS) quantities of the semi-Gaussian beam in Figure (\ref{fig:rms}). Although the RMS quantities of the semi-Gaussian
beam are oscillatory, they remain close to the quantities of the associated K-V beam and they converge as the resolution increases. 
As mentioned by the other authors \cite{Crouseilles08II}, the oscillatory behavior is due to the fact that the semi-Gaussian is not exactly a steady state distribution.

\section{Conclusion}
In this paper, we have presented the adaptive remapped PIC method to the high-dimensional Vlasov equation and demonstrated in
linear Landau damping, the two stream instability, and the beam propagation problems. The new method reduces the numerical noise 
significantly. There are two extensions of the current research. The first will be the development of a scalable
algorithm based on domain decomposition in phase space. The second will be the introduction of time adaptivity to the current algorithm 
that the hierarchy of locally-refined grids are dynamically created from the particle distribution at every remapping step.

\section*{Acknowledgments}
This work was supported by the U.S. Department of Energy Office of
Advanced Scientific Computing Research under contract number DE-AC02-05CH11231 at the Lawrence Berkeley National Laboratory.
G.~H.~Miller was supported by DOE contract number DE-SC0001981.

\bibliographystyle{siam}
\bibliography{vlasovpic2d}

\end{document}